\documentclass[12pt]{article}
\usepackage{amsmath,amssymb,amscd,amsthm,bm}
\usepackage{amsfonts}
\usepackage{graphicx}
\newtheorem{theorem}{Theorem}[section]
\newtheorem{proposition}[theorem]{Proposition}

\newtheorem{lemma}[theorem]{Lemma}
\newtheorem{remark}[theorem]{Remark}

\numberwithin{equation}{section}
\begin{document}

\begin{center}
\Large
Quadratic Embedding Constants of Block Graphs
and Characterization of Cactus-like Graphs
\end{center}

\bigskip

\begin{center}
Wojciech M{\l}otkowski \\
Instytut Matematyczny \\
Uniwersytet Wroc{\l}awski \\
Plac Grunwaldzki 2/4, 50-384 Wroc{\l}aw, Poland \\
mlotkow@math.uni.wroc.pl
\\
and
\\
Nobuaki Obata\\
Center for Data-driven Science and Artificial Intelligence \\
Tohoku University\\
Sendai 980-8576 Japan \\
obata@tohoku.ac.jp
\end{center}

\bigskip

\begin{quote}
\textbf{Abstract}\enspace
A (finite or infinite) connected graph is called \textit{cactus-like} 
if every pair of distinct maximal cliques has
at most one vertex in common and its clique graph is a tree.
Equivalently, a cactus-like graph is a block graph whose clique graph is a tree.
We prove that a finite connected graph $G$ with at least two vertices
is cactus-like if and only if its quadratic embedding constant satisfies
$\mathrm{QEC}(G)<-1/2$.
Moreover, we establish upper and lower bounds for $\mathrm{QEC}(G)$
in terms of the clique degree of a block graph $G$.
For the proofs, we develop a method for constructing a conditionally
negative definite kernel on a block graph by concatenating local kernels
along geodesics in the block graph.
Concrete examples include free products of complete graphs and trees.
\end{quote}

\begin{quote}
\end{quote}

\begin{quote}
\end{quote}

\section{Introduction}

For a (finite or infinite) connected graph $G=(V,E)$ with $|V|\ge2$ the 
\textit{quadratic embedding constant} (QEC for short), denoted by $\mathrm{QEC}(G)$, 
is defined to be the conditional supremum of the canonical quadratic
function: 
\[
\langle f,Df \rangle=\sum_{x,y\in V} d_G(x,y)f(x)f(y),
\qquad f\in C_0(V),
\]
associated with the distance matrix
$D=[d_G(x,y)]$ of $G$ subject to
\[
\langle f,f \rangle=1, \quad \langle \bm{1},f \rangle=0,
\quad f\in C_0(V),
\]
where $C_0(V)$ is the space of real
column vectors $f=[f(x)]_{x\in V}$ with finite supports,
$\bm{1}\in C(V)$ is the all-one column vector,
and $\langle\cdot,\cdot\rangle$ the canonical inner product.
If $G$ is finite, the conditional supremum is attained at some
$f\in C(V)=C_0(V)$ so that we may write
\begin{equation}\label{01eqn:def of QEC for a finite graph}
\mathrm{QEC}(G)
=\max\{\langle f,Df \rangle\,;\, f\in C(V), \,
\langle f,f \rangle=1, \, \langle \bm{1},f \rangle=0\}.
\end{equation}
The QEC becomes a numeric invariant of connected graphs
on two or more vertices,
and takes values in $[-1,+\infty)$.

Since the QEC was introduced in \cite{Obata-Zakiyyah2018},
research interest in it has steadily increased,
with most studies focusing on finite graphs;
see \cite{Obata2026-book} for a comprehensive survey of recent results.
Concrete values of QEC are known for some special classes of graphs,
e.g., paths \cite{Mlotkowski2022},
fan graphs \cite{Mlotkowski-Obata2025a,Mlotkowski-Obata2025b},
theta graphs \cite{Mlotkowski-Skrzypczyk-Wojtylak2026,Skrzypczyk2026},
wheels \cite{Obata2017} and others \cite{Irawan-Sugeng2021,
Purwaningsih-Sugeng2021}.
Formulas and estimates for QEC under various graph operations
have also been investigated, 
e.g., for cluster graphs \cite{Choudhury-Nandi2024},
corona graphs \cite{Ferdi-Baskoro-Obata-Santika2026},
graph joins \cite{Lou-Obata-Huang2022},
star products \cite{Mlotkowski-Obata2020},
lexicographic products \cite{Mou2026}, and so on.
The investigation of graphs with a positive QEC has only recently begun
\cite{Mlotkowski-Skrzypczyk-Wojtylak2026,Obata2023a,Obata2023b}.
The classification of graphs in terms of QEC
remains a central open problem, 
though progress has been limited so far,
see \cite{Baskoro-Obata2021,Baskoro-Obata2024}.
In this paper we present a result
that is expected to constitute an important step toward this goal.

The most fundamental fact that motivated 
the introduction of the QEC comes from the famous works of Schoenberg 
\cite{Schoenberg1935,Schoenberg1938}.
Let $G=(V,E)$ be a (finite or infinite) connected graph and
$\mathcal{H}$ a Hilbert space.
A map $\varphi:V\rightarrow \mathcal{H}$ is called
a \textit{quadratic embedding} of $G$ if
\[
\|\varphi(x)-\varphi(y)\|^2=d_G(x,y),
\qquad x,y\in V,
\]
where $\|\cdot\|$ stands for the norm of $\mathcal{H}$.
The concept of quadratic embedding has been developed
into a central one along with the Euclidean distance geometry
\cite{Alfakih2018,Balaji-Bapat2007,Jaklic-Modic2013,Jaklic-Modic2014,
Liberti-Lavor-Maculan-Mucherino2014}.
In our terminology we recall the following

\begin{theorem}[Schoenberg]
\label{01thm:Schoenberg}
For a (finite or infinite) connected graph
$G=(V,E)$ with $|V|\ge2$ the following
three conditions are equivalent:
\begin{enumerate}
\item[\upshape (i)] $G$ admits a quadratic embedding;
\item[\upshape (ii)] the distance matrix $D$ of $G$ is
conditionally negative definite, i.e., 
$\langle f,Df\rangle\le 0$ for any $f\in C_0(V)$
with $\langle \bm{1},f\rangle=0$;
\item[\upshape (iii)] $\mathrm{QEC}(G)\le 0$.
\end{enumerate}
\end{theorem}

The essential point lies in the equivalence 
(i) $\Leftrightarrow$ (ii),
which was first proved in the seminal works of Schoenberg
\cite{Schoenberg1935,Schoenberg1938},
while the equivalence (ii) $\Leftrightarrow$ (iii) is
verified easily from the definition of the QEC.
As a result, the QEC allows us to determine numerically
whether a connected graph admits a quadratic embedding.
In this context, we will determine the class of graphs with
$\mathrm{QEC}(G)<-1/2$,
which represents the first non-trivial achievement
after Schoenberg's result (Theorem~\ref{01thm:Schoenberg}).
We will prove the following

\begin{theorem}\label{01thm:main result}
Let $G=(V,E)$ be a finite connected graph with $|V|\ge2$.
Then the following conditions are equivalent:
\begin{enumerate}
\item[\upshape (i)] $G$ is cactus-like,
that is, every pair of distinct maximal cliques has 
at most one vertex in common and the clique graph $\Gamma(G)$ is a tree;
\item[\upshape (ii)]  $G$ is a block graph
whose clique graph $\Gamma(G)$ is a tree;
\item[\upshape (iii)] $\mathrm{QEC}(G)<-1/2$.
\end{enumerate}
\end{theorem}

\begin{figure}[hbt]
\begin{center}
\includegraphics[width=360pt]{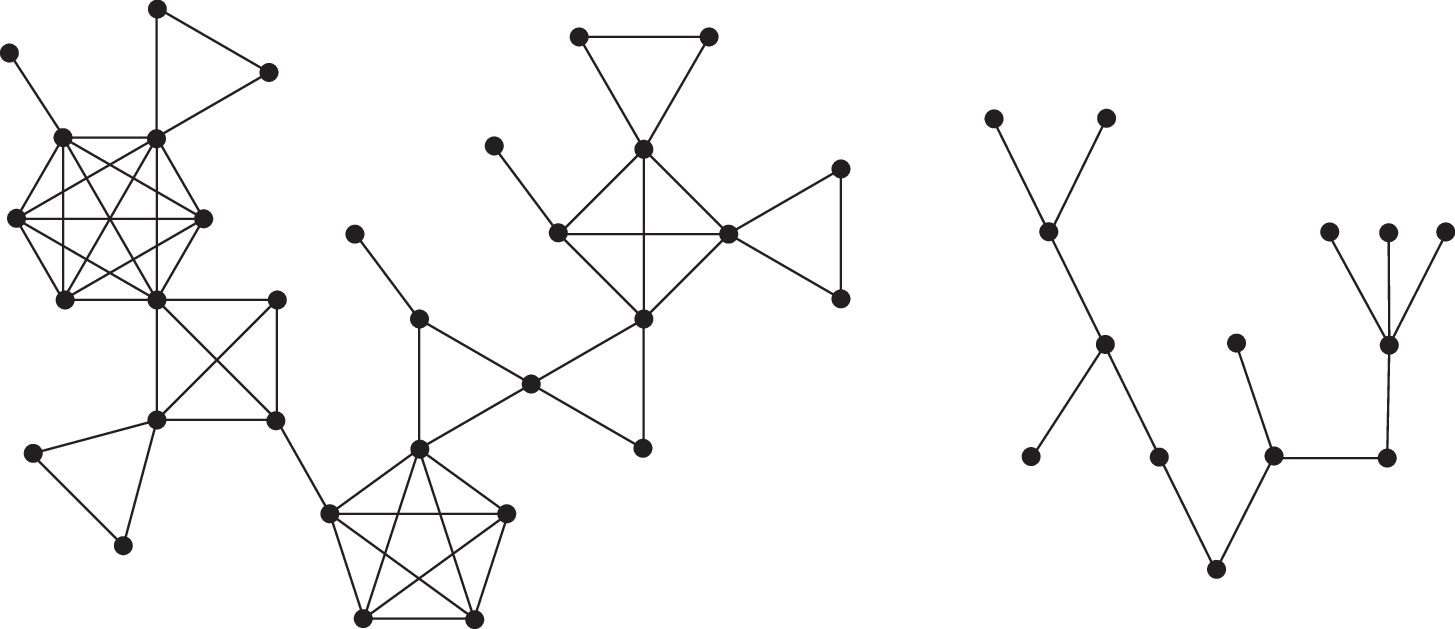}
\caption{Cactus-like graph $G$ and its clique graph $\Gamma(G)$}
\label{fig:Cactus-like}
\end{center}
\end{figure}

A cactus-like graph was introduced 
in the recent paper \cite{Baskoro-Obata2024},
where the implication (iii) $\Rightarrow$ (i) was proved.
The main contribution of this paper is 
the converse (i) $\Rightarrow$ (iii), which has been left open.
Meanwhile, the equivalence (i) $\Leftrightarrow$ (ii),
which can be verified directly by definition, provides a clear structural
characterization of cactus-like graphs.
To gain an intuitive understanding of a cactus-like graph,
an example is given in Fig.\ref{fig:Cactus-like}.

In the proof of the implication (i) $\Rightarrow$ (iii), 
two results play a key role.
One is a generalization of a result on
the length function of trees or of free product groups
studied within harmonic analysis
\cite{Mlotkowski1988,Mlotkowski1989,Mlotkowski1992,Valette1992}.
The other is a technical estimate of QEC of star products 
obtained in \cite{Baskoro-Obata2021}.
Thus, our study of QEC of block graphs led us to
a link of harmonic analysis on free products.
Moreover, determining the infinite connected graphs
that satisfy $\mathrm{QEC}(G)<-1/2$ is an interesting open question and
the subject of ongoing research.

This paper is organized as follows.
In Section 2 we recall some notions and notations of graphs
in relation to QEC, and assemble some known results necessary for
our goal.
In Section 3 we describe some properties of block graphs
and introduce cactus-like graphs in this context.
In particular, we obtain equivalent conditions for a graph
to be cactus-like
(Proposition~\ref{03prop:equivalent conditions for cactus-like}),
from which we see that two conditions
(i) and (ii) in Theorem \ref{01thm:main result} are equivalent.
In Section 4 we derive upper and lower bounds for QEC of a block graph
and complete the proof of Theorem~\ref{01thm:main result}.
The upper bound of QEC of a finite block graph
with two or more maximal cliques is also noteworthy 
(Theorem~\ref{04thm:upper bound for finite block graph}).
In Section 5 we prove the key result on 
conditionally negative definite matrices on a block graph
(Theorem~\ref{04thm:psi is CND}),
which plays an essential role in the proof of the main result.
The proof proceeds by constructing
a conditionally negative definite kernel by
concatenating local ones along geodesics in a block graph.
This method is an extension of a similar idea established
for trees and Cayley graphs of free product groups
\cite{Mlotkowski1992}.
In Section 6 we discuss trees and free products of complete graphs
as special cases and make some comments.
In Appendix we review the formula for QEC of a 
star product of complete graphs,
which is used for the lower bound for QEC of a block graph.

\section{Preliminaries}
\label{02sec:Preliminaries}

\subsection{Graphs and Subgraphs}

A \textit{graph} is a pair $G=(V,E)$ consisting of a non-empty set $V$ and
a subset of the two-element subsets of $V$.
An element of $V$ is called a \textit{vertex}
and an element of $E$ an \textit{edge}.
A graph $G=(V,E)$ is called \textit{finite} if $V$ is a 
(non-empty) finite set,
and \textit{infinite} otherwise.
In this paper we deal with both finite and infinite graphs.

Let $G=(V,E)$ be a graph.
Two vertices $x,y\in V$ are called \textit{adjacent}
if $\{x,y\}\in E$.
In that case, we write $x\sim y$ as usual.
Note that this definition of a graph precludes loops,
that is, $x\not\sim x$.
For any $x,y\in V$
a finite sequence of vertices $x_0,x_1,\dots,x_m\in V$ such that
\begin{equation}\label{01eqn:a walk}
x=x_0\sim x_1\sim\dots\sim x_m=y\,,
\end{equation}
where $m\ge0$ and some of $x_0,x_1,\dots,x_m$ may coincide,
is called a \textit{walk of length $m$} 
connecting $x$ and $y$.
If $x_0, x_1, \dots, x_m$ are mutually distinct,
the walk in \eqref{01eqn:a walk} is called a \textit{path of length $m$}.
If $m\ge3$, $x=y$, and $x_0,x_1,\dots,x_{m-1}$ are mutually distinct,
then the walk in \eqref{01eqn:a walk} is called
a \textit{cycle of length $m$}.

A graph $G=(V,E)$ is called \textit{connected} if
any two vertices $x,y\in V$ are connected by a walk.
In that case, the length of a shortest walk 
(equivalently, a shortest path) connecting $x$ and $y$
is called the \textit{graph distance} between $x$ and $y$,
and is denoted by $d_G(x,y)$ or simply
by $d(x,y)$ when there is no risk of confusion.

Let $G=(V,E)$ and $H=(W,F)$ be two graphs.
We say that $H$ is a \textit{subgraph} of $G$
and write $H\subset G$ if $W\subset V$ and $F\subset E$.
A subgraph $H=(W,F)$ of $G=(V,E)$ is called an
\textit{induced subgraph} if $F=E\!\restriction_W$ holds, 
that is, for any $x,y\in W$ 
we have $\{x,y\}\in F \,\Leftrightarrow\,\{x,y\}\in E$.
In that case we also say that $H$ is the subgraph induced by $W$,
and write $H=\langle W\rangle$.

Let $G=(V,E)$ be a connected graph and
$H=(W,F)$ a connected subgraph of $G$.
Let $d_G$ and $d_H$ be the graph distances of $G$ and $H$,
respectively.
We say that $H$ is an \textit{isometric subgraph} of $G$ if
\begin{equation}\label{03eqn:isometric subgraphs}
d_H(x,y)=d_G(x,y),
\qquad x,y\in W.
\end{equation}
An isometric subgraph is necessarily an induced subgraph,
but the converse is false.

\subsection{Quadratic Embedding Constant}

Let $G=(V,E)$ be a connected graph with $|V|\ge2$
and $d(x,y)$ the graph distance between two vertices $x$ and $y$
in $V$.
Recall that $G$ may be finite or infinite.
Let $C(V)$ be the space of real column vectors $f=[f(x)]_{x\in V}$
and $C_0(V)$ its subspace of those with finite supports,
that is, those $f=[f(x)]$ with $f(x)=0$ except finitely many $x\in V$.

The canonical quadratic function associated to
the distance matrix $D=[d(x,y)]_{x,y\in V}$ of $G$ is defined by
\begin{equation}\label{02eqn:canonical quadratic function ass to D}
\langle f,Df \rangle=\sum_{x,y\in V} d(x,y)f(x)f(y),
\qquad f\in C_0(V),
\end{equation}
where $\langle\cdot,\cdot\rangle$ stands for the
canonical inner product on $C_0(V)$.
The \textit{quadratic embedding constant} 
(\textit{QEC} for short) of $G$ is 
the conditional supremum of $\langle f,Df \rangle$ defined by
\begin{equation}\label{02eqn:def of QEC}
\mathrm{QEC}(G)
=\sup\{\langle f,Df \rangle\,;\, f\in C_0(V), \,
\langle f,f \rangle=1, \, \langle \bm{1},f \rangle=0\},
\end{equation}
where $\bm{1}\in C(V)$ is the column vector whose entries are all one
and $\langle \bm{1},f \rangle$ is understood as a natural extension
of the inner product on $C_0(V)$.
If $G=(V,E)$ is a finite connected graph with $|V|=n\ge2$,
we have $C(V)=C_0(V)\cong \mathbb{R}^n$ and
$\mathrm{QEC}(G)$ is attained at some $f\in C(V)$
satisfying $\langle f,f \rangle=1$ and $\langle \bm{1},f \rangle=0$.
Hence, \eqref{02eqn:def of QEC} becomes
\begin{equation}\label{02eqn:def of QEC for a finite graph}
\mathrm{QEC}(G)
=\max\{\langle f,Df \rangle\,;\, f\in C(V), \,
\langle f,f \rangle=1, \, \langle \bm{1},f \rangle=0\},
\end{equation}
as stated in Introduction.
The QEC is a numeric invariant of connected graphs
on two or more vertices.

The following result, the proof of which is straightforward from definition,
will be used freely without further comment.

\begin{proposition}[\cite{Obata-Zakiyyah2018}]
\label{03prop:isometric subgraphs}
Let $G=(V,E)$ be a connected graph 
and $H=(W,F)$ its isometric subgraph with $|W|\ge2$.
Then we have $\mathrm{QEC}(H)\le \mathrm{QEC}(G)$.
\end{proposition}

The set of integers $\mathbb{Z}=\{\dots, -2,-1,0,1,2,\dots\}$
equipped with the natural adjacency relation $i\sim i+1$
becomes an (infinite) graph called
the \textit{two-sided infinite path}
and is denoted by $\mathbb{Z}$ for simplicity.
The subgraph induced by $\{0,1,2,\dots\}$ is
called the \textit{one-sided infinite path}
and denoted by $\mathbb{Z}_+$.

\begin{proposition}[\cite{Mlotkowski2022,Mlotkowski-Obata2020}]
\label{03prop:QEC(Z)}
It holds that 
\[
\lim_{n\rightarrow\infty} \mathrm{QEC}(P_n)
=\mathrm{QEC}(\mathbb{Z}_+)
=\mathrm{QEC}(\mathbb{Z})=-\frac12,
\]
where $P_n$ is a path graph on $n$ vertices.
\end{proposition}

The original proof \cite{Mlotkowski-Obata2020}
requires a precise estimate of 
$\mathrm{QEC}(P_n)$ without knowing its explicit value.
Subsequently, the explicit value
\[
\mathrm{QEC}(P_n)
=-\frac{1}{1+\cos \pi/n},
\qquad
n=2,3,\dots,
\]
was determined through a skillful analysis of a certain determinant \cite{Mlotkowski2022}.

Following general notation of graph theory, a finite graph
on $n\ge1$ vertices is called \textit{complete} if
every pair of distinct vertices is adjacent,
and such a graph is denoted by $K_n$ as usual.
To avoid ambiguity, we always use the term `complete graph'
exclusively for finite graphs,
and refer to an infinite graph 
in which every pair of distinct vertices is adjacent
as an `infinite complete graph' to distinguish it from 
the former.
An infinite complete graph is denoted by $K_\infty$,
regardless of the cardinality of the vertex set.
The distinction between infinite cardinalities plays 
no role throughout this paper.

\begin{proposition}\label{02prop:QEC(K_n)}
It holds that 
\[
\mathrm{QEC}(K_n)
=\mathrm{QEC}(K_\infty)=-1,
\qquad
n\ge2.
\]
\end{proposition}

In fact, for a finite or infinite complete graph $G=(V,E)$
with distance matrix $D$ we have
$\langle f,Df\rangle=-1$ for any
$f\in C_0(V)$ with $\langle f,f\rangle=1$
and $\langle \bm{1},f\rangle=0$.

\subsection{Star Products}

A graph $G=(V,E)$ equipped with a distinguished vertex $o\in V$
is called a \textit{rooted graph} with a root $o$
and is denoted by $(G,o)=(V,E,o)$.
For two rooted graphs $(G_1,o_1)=(V_1,E_1,o_1)$
and $(G_2,o_2)=(V_2,E_2,o_2)$ with $V_1\cap V_2=\emptyset$
the graph obtained by gluing their roots is called 
the \textit{star product} of $(G_1,o_1)$ and $(G_2,o_2)$,
and is denoted by $(G_1,o_1)\star (G_2,o_2)$.
If there is no danger of confusion, 
we omit reference to the roots $o_1$ and $o_2$.

\begin{proposition}[\cite{Obata-Zakiyyah2018}]
\label{03prop:star product}
Let $G_1$ and $G_2$ be connected rooted graphs
on two or more vertices,
and let $G_1\star G_2$ be their star product.
Then we have
\[
\max\{\mathrm{QEC}(G_1), \mathrm{QEC}(G_2)\}
\le \mathrm{QEC}(G_1\star G_2).
\]
Moreover, 
if $\mathrm{QEC}(G_1)\le 0$ and $\mathrm{QEC}(G_2)\le 0$, we have
\[
\mathrm{QEC}(G_1\star G_2)\le 0.
\]
\end{proposition}

The first half is immediate from
Proposition~\ref{03prop:isometric subgraphs}
since both $G_1$ and $G_2$ are isometric subgraphs of $G_1\star G_2$.
The second half is proved by constructing 
a quadratic embedding map of $G_1\star G_2$,
see also \cite{Mlotkowski-Obata2020}
for a computational proof.
It is an interesting problem
to obtain precise estimate of $\mathrm{QEC}(G_1\star G_2)$ 
for two connected graphs $G_1$ and $G_2$ with
$\mathrm{QEC}(G_1)< 0$ and $\mathrm{QEC}(G_2)< 0$.
In this connection the following result is useful to our goal.

\begin{proposition}[\cite{Baskoro-Obata2021}]
\label{02prop:Baskoro-Obata2021}
Let $G=(V,E)$ be a finite connected graph with $|V|\ge2$
with $\mathrm{QEC}(G)<0$.
Let $o\in V$.
If there exists $f\in C(V)$ satisfying
\[
\mathrm{QEC}(G)=\langle f, Df\rangle,
\quad
\langle f, f\rangle=1,
\quad
\langle \bm{1}, f\rangle=0,
\quad
f(o)\neq0,
\]
then for any $k\ge2$ we have
\[
\mathrm{QEC}(G)<\mathrm{QEC}((G,o)\star K_k)<0,
\]
where $K_k$ is the complete graph on $k$ vertices
and its root is chosen arbitrarily.
\end{proposition}

\section{Block Graphs}
\label{03sec:Block Graphs}

Our terminology concerning block graphs follows standard textbooks
\cite{Balakrishnan-Ranganathan2012,Bondy-Murty2008},
see also \cite{Harary1963} for the original idea of Harary.
Although many references restrict their attention to finite graphs, 
this paper also considers infinite graphs.
To avoid any ambiguity arising from this broader setting, 
we provide some explanations,
even at the expense of some redundancy.
Recall that throughout this paper, 
the term `graph' refers to either a finite or an infinite graph.

\subsection{2-Connected Components}

A connected graph $G=(V,E)$ is called \textit{2-connected}
if the induced subgraph $G\setminus\{x\}=\langle V\setminus\{x\}\rangle$ remains
connected for any $x\in V$.
By definition, the singleton graph $K_1$ is connected
but not 2-connected,
while the complete graph $K_n$ on $n\ge2$ vertices is 2-connected.
A 2-connected subgraph $H=(W,F)$ of a graph $G$ is called
a \textit{2-connected component} if
there is no 2-connected subgraph $H_1$ of $G$ 
such that $H\subset H_1$ and $H\neq H_1$.
A 2-connected component is necessarily an induced subgraph.

The following assertions follow directly from
the definitions via elementary arguments.

\begin{lemma}\label{03lem:intersecting 2-connected components}
Let $G=(V,E)$ be a graph.
Let $H_1=(W_1,F_1)$ and $H_2=(W_2,F_2)$ be
2-connected subgraphs of $G$.
If $|W_1\cap W_2|\ge2$, then
$H_1\cup H_2=(W_1\cup W_2,F_1\cup F_2)$ is also a
2-connected subgraph of $G$.
\end{lemma}

\begin{lemma}\label{03lem:2-connected is contained in unique component}
Let $G=(V,E)$ be a graph.
Any 2-connected subgraph $H_0$ of $G$
is contained in a unique 2-connected component.
In particular, any edge of $G$ 
is contained in a unique 2-connected component.
\end{lemma}

\begin{lemma}\label{03lem:H1 cap H2 <=1 for 2-connected component}
Let $G=(V,E)$ be a graph.
For any distinct 2-connected components
$H_1=(W_1,F_1)$ and $H_2=(W_2,F_2)$ we have
$|W_1\cap W_2|\le1$.
\end{lemma}

\subsection{Maximal Cliques and Clique Graphs}

Let $G=(V,E)$ be a graph.
A non-empty (finite or infinite) subset $C\subset V$ is
called a \textit{clique} of $G$ if every pair of distinct vertices
in $C$ is adjacent,
or equivalently if $\langle C \rangle$ is a 
finite or infinite complete graph.
By a slight abuse of notation, 
the term `clique' may also refer to the subgraph it induces.

\begin{lemma}\label{01lem:clique in 2-connected component}
Let $G=(V,E)$ be a graph.
Let $H=(W,F)$ be a 2-connected component and
$C\subset V$ a clique.
If $|W\cap C|\ge2$, we have $C\subset W$.
\end{lemma}

\begin{proof}
From $|W\cap C|\ge2$ we see that $\langle C\rangle$ is
a finite or infinite complete graph on two or more vertices,
hence is 2-connected.
It then follows by
Lemma~\ref{03lem:intersecting 2-connected components} that
$H\cup \langle C\rangle$ is a 2-connected subgraph of $G$.
But since $H$ is a maximal 2-connected subgraph by assumption,
we have $H=H\cup \langle C\rangle$,
which implies that $W=W\cup C$, that is, $C\subset W$.
\end{proof}

A clique of a graph $G=(V,E)$ is called \textit{maximal} if
it is maximal in the family of cliques with respect to the
set inclusion relation.
The following assertion is obvious for a finite graph $G$;
however, for a general infinite graph $G$,
we need to invoke the Kuratowski--Zorn lemma.

\begin{lemma}\label{03lem:existence of maximal clique}
Let $G=(V,E)$ be a graph.
For any clique $C_0$ of $G$ there exists a maximal clique
$C$ such that $C_0\subset C$.
In particular,
any vertex as well as any edge is contained in at least one 
maximal clique.
\end{lemma}

The set $\mathcal{C}$ of
all maximal cliques of a graph $G$ induces a graph,
where two distinct maximal cliques $C_1,C_2\in\mathcal{C}$ are 
defined to be adjacent if $C_1\cap C_2 \neq\emptyset$.
The graph thus obtained is called the \textit{clique graph} of $G$,
and is denoted by $\Gamma(G)$.

\begin{lemma}\label{01lem:Gamma(G) is connected}
Let $G=(V,E)$ be a graph.
Then the clique graph $\Gamma(G)$ is connected 
if and only if so is $G$.
\end{lemma}

The proof follows directly
from the observation that any edge of $G$ is contained
in a maximal clique, see Lemma~\ref{03lem:existence of maximal clique}.
For later convenience, we mention some detailed results
on paths in $G$ and in $\Gamma(G)$.

\begin{lemma}\label{03lem:corresponding paths (1)}
Let $G=(V,E)$ be a graph and $m\ge1$.
Let $x_0\sim x_1\sim\dotsb\sim x_m$ be
a shortest path in $G$ connecting $x_0$ and $x_m$.
For $i=1,2,\dots,m$ let $C_i$ be a maximal clique containing
$\{x_{i-1},x_i\}$.
Then $C_1\sim C_2\sim \dots \sim C_m$ is a shortest path
in $\Gamma(G)$ connecting $C_1$ and $C_m$.
\end{lemma}

\begin{lemma}\label{03lem:corresponding paths (2)}
Let $G=(V,E)$ be a graph and $m\ge1$.
Let $C_0\sim C_1\sim \dots \sim C_m$ be a shortest path
connecting $C_0$ and $C_m$ in the clique graph $\Gamma(G)$.
For $i=1,2,\dots,m$ take a vertex $x_i\in C_{i-1}\cap C_i$.
Then $x_1\sim x_2\sim \dots \sim x_m$ is a shortest path
in $G$ connecting $x_1$ and $x_m$.
\end{lemma}

The proofs of the above assertions require
only elementary combinatorial arguments directly from
the definition and are therefore omitted,
for relevant results see \cite{Dutton-Brigham1995},
and for more on clique graphs see \cite{Hamelink1968,
Howorka1979,Roberts-Spencer1971,Szwarcfiter2003}.

\subsection{Block Graphs}

A connected graph $G=(V,E)$ is called a \textit{block graph}
if every 2-connected component is a clique.
The singleton graph $K_1$ is a block graph by definition.
Any tree is a block graph, where any 2-connected component
is an edge and vice versa.

\begin{lemma}\label{01lem:2-connected comp = max clique}
Let $G=(V,E)$ be a block graph with $|V|\ge2$.
Then a subset $W\subset V$ is a maximal clique
if and only if $\langle W\rangle$ is a 2-connected component.
\end{lemma}

\begin{proof}
Suppose first that $\langle W \rangle$ is a 2-connected component.
By definition of a block graph, $W$ is a clique.
Let $W_1$ be a clique of $G$ such that $W\subset W_1$.
Since $\langle W\rangle$ is a 2-connected graph, we have $|W|\ge2$.
Then $\langle W_1\rangle$ is a finite or infinite complete graph
on $|W_1|\ge2$ vertices, and hence is 2-connected.
Since $\langle W \rangle$ is a maximal 2-connected subgraph,
we come to $W=W_1$, which means that
$W$ is a maximal clique.

Conversely, let $W\subset V$ be a maximal clique.
Since $V$ contains no isolated vertices by $|V|\ge2$,
$\langle W\rangle$ is a finite or infinite 
complete graph on $|W|\ge2$ vertices and is therefore 2-connected.
By Lemma~\ref{03lem:2-connected is contained in unique component}
there exists a 2-connected component $H_1=(W_1,F_1)$ 
containing $\langle W\rangle$, that is, $W\subset W_1$.
Note that $W_1$ is a clique by definition of a block graph.
Since $W$ is a maximal clique, we conclude that $W=W_1$ and 
$\langle W\rangle$ is a 2-connected component.
\end{proof}

\begin{lemma}\label{01lem:H1 cap H2 <=1 for block graph}
Let $G=(V,E)$ be a block graph.
Then for any distinct maximal cliques $C_1$ and $C_2$ we have
$|C_1\cap C_2|\le1$.
\end{lemma}

\begin{proof}
Immediate from
Lemmas \ref{03lem:H1 cap H2 <=1 for 2-connected component}
and \ref{01lem:2-connected comp = max clique}.
\end{proof}

\begin{remark}\normalfont
The converse assertion in Lemma~\ref{01lem:H1 cap H2 <=1 for block graph}
does not hold.
Consider a cycle $G$ on four or more vertices.
Then $G$ itself is 2-connected, hence is a unique 2-connected component.
But this 2-connected component is not a clique,
which means that $G$ is not a block graph.
On the other hand, since each edge of $G$ is a maximal clique,
we have $|C_1\cap C_2|\le1$ for any distinct maximal cliques of $G$.
\end{remark}

\begin{proposition}\label{03prop:Block graph is geodetic}
Any block graph is geodetic,
i.e., every pair of vertices is
connected by a unique shortest path (called a geodesic).
\end{proposition}

\begin{proof}
The proof can be given by an elementary argument based on the definition.
For example, it can be proved by induction on $n\ge1$ that
any pair of vertices $x,y$ with $d(x,y)=n$ admits
a unique shortest path.
\end{proof}

The result stated in Proposition~\ref{03prop:Block graph is geodetic}
is well-known and can be found, 
for example, in \cite{Behtoei-Jannesari-Taeri2010,Dobrynin-Gutman1994}.
However, their main purpose was to obtain
an interesting characterization of the block graphs
using the Wiener and Szeged indices.

Let $G=(V,E)$ be a graph. 
The \textit{clique degree} of a vertex $x\in V$
is defined to be the number of maximal cliques containing $x$,
and is denoted by $\mathrm{cd}(x)$. 
Obviously, $\mathrm{cd}(x)\ge1$ for any $x\in V$.
If there are infinitely many maximal cliques containing $x$,
we simply write $\mathrm{cd}(x)=\infty$.
The clique degree is useful for the structure of block graphs.

\subsection{Cactus-like Graphs}
\label{03subsec:Cactus-like Graphs}

Although the concept of a cactus-like graph introduced in \cite{Baskoro-Obata2024}
is limited to finite graphs,
the definition applies equally well to infinite graphs
without modification.
A graph $G=(V,E)$ is called \textit{cactus-like} if
the clique graph $\Gamma(G)$ is a tree and 
$|C_1\cap C_2|\le1$ for any pair of distinct maximal cliques
$C_1$ and $C_2$ of $G$.
Any cactus-like graph is connected
by Lemma~\ref{01lem:Gamma(G) is connected}.
By definition, the complete graph $K_n$ with $n\ge1$
as well as infinite complete graphs $K_\infty$ is cactus-like.

\begin{proposition}\label{03prop:cactus-like graph is a block graph}
Any cactus-like graph is a block graph.
\end{proposition}

\begin{proof}
Let $G=(V,E)$ be a cactus-like graph
and $H=(W,F)$ any 2-connected component of $G$.
Since $H$ is an induced subgraph, the adjacency relation
in $H$ coincides with that in $G$.
We need to show that $W$ is a clique.
If $|W|=2$, then $H\cong K_2$ and hence $W$ is certainly a clique.
Hereafter we assume that $|W|\ge3$.
It is sufficient to prove that
every pair of distinct vertices in $W$ is adjacent.
We prove this by contradiction.

Suppose that
$H=(W,F)$ is a 2-connected component of $G$
with $|W|\ge3$ that is not a clique.
Then there exists an induced path of length 2 in $H$, say, 
\begin{equation}\label{03eqn:in proof Prop 3.13 path of length 2}
a \sim c \sim b,
\end{equation}
where $a,b,c\in W$ are mutually distinct and $a\not\sim b$.
Let $C_1$ and $C_2$ be maximal cliques of $G$ containing
$\{a,c\}$ and $\{c,b\}$, respectively.
Such maximal cliques
exist by Lemma~\ref{03lem:existence of maximal clique}.
Since \eqref{03eqn:in proof Prop 3.13 path of length 2}
becomes a shortest path connecting $a$ and $b$ in $G$,
by Lemma~\ref{03lem:corresponding paths (1)} we see that
\begin{equation}\label{03eqn:C-path in proof Prop 3.12}
C_1\sim C_2
\end{equation}
is a shortest path, in particular,
$C_1\neq C_2$ and $C_1\cap C_2=\{c\}$.

Since $H\setminus\{c\}$ is connected
because $H$ is 2-connected, 
we may take a shortest path connecting $a$ and $b$ in $H\setminus\{c\}$,
say,
\begin{equation}\label{03eqn:in proof Prop 3.12 (1)}
a=x_0\sim x_1\sim \dotsb\sim x_{s-1}\sim x_s=b,
\qquad
x_i\in W\setminus\{c\},
\quad s\ge2.
\end{equation}
For $i=1,2,\dots,s$ let $D_i$ be a maximal clique containing
$\{x_{i-1},x_i\}$.
It is easy to see that
\begin{equation}\label{03eqn:D-path in proof Prop 3.12}
D_1\sim D_2\sim\dotsb \sim D_s
\end{equation}
and these are mutually distinct.
In fact, if $D_i=D_j$ happens for $1\le i<j\le s$,
we have $\{x_{i-1},x_j\}\subset D_i=D_j$ and $x_{i-1}\sim x_j$.
This contradicts that \eqref{03eqn:in proof Prop 3.12 (1)}
is a shortest path connecting $a$ and $b$ in $H\setminus\{c\}$.

In view of $a\in C_1\cap D_1$ and $b\in C_2\cap D_s$
there are three cases:
\par
(Case 1) $C_1=D_1$ and $C_2=D_s$.
Note first that $s\ge3$.
In fact, if $s=2$ we have $\{c,x_1\}\subset C_1=D_1$
and $\{c,x_1\}\subset C_2=D_2$.
Then we come to $|C_1\cap C_2|\ge2$, which contradicts that
$G$ is cactus-like.
Concatenating the two paths
in \eqref{03eqn:C-path in proof Prop 3.12}
and \eqref{03eqn:D-path in proof Prop 3.12},
we obtain a cycle in $\Gamma(G)$ of length $s\ge3$.
But this contradicts that $\Gamma(G)$ is a tree.

\par
(Case 2) $C_1\sim D_1$ and $C_2=D_s$.
We note that $D_i\neq C_1$ for $2\le i\le s$.
In fact, if $D_i= C_1$ occurs, we have
$\{x_0,x_i\}\subset D_i= C_1$ and hence $x_0\sim x_i$.
This contradicts that \eqref{03eqn:in proof Prop 3.12 (1)} 
is a shortest path connecting $a$ and $b$ in $H\setminus\{c\}$.
Thus, we come to
\[
C_1\sim D_1 \sim D_2 \sim \dots \sim D_s=C_2
\]
and these are mutually distinct.
Then combining \eqref{03eqn:C-path in proof Prop 3.12},
we obtain a cycle of length $s+1\ge3$ in $\Gamma(G)$.
But this contradicts that $\Gamma(G)$ is a tree.
Similarly, the case where $C_1= D_1$ and $C_2\sim D_s$
leads to a contradiction.

\par
(Case 3) $C_1\sim D_1$ and $C_2\sim D_s$.
In a similar argument as in (Case 2) we obtain 
a cycle of length $s+2\ge4$ in $\Gamma(G)$,
which is a contradiction.

\par
Consequently, in any case we come to a contradiction
and the proof is complete.
\end{proof}

We are in a good position to
collect the equivalent conditions
for a graph to be cactus-like.
In particular,
the equivalence of two conditions
(i) $\Leftrightarrow$ (ii) in the main Theorem
\ref{01thm:main result} is established.

\begin{proposition}\label{03prop:equivalent conditions for cactus-like}
For a graph $G=(V,E)$ the following conditions are equivalent:
\begin{enumerate}
\item[\upshape (i)] $G$ is cactus-like,
that is, every pair of distinct maximal cliques has 
at most one vertex in common and the clique graph $\Gamma(G)$ is a tree;
\item[\upshape (ii)] $G$ is a block graph and $\Gamma(G)$ is a tree.
\item[\upshape (iii)] $G$ is a block graph and
$\mathrm{cd}(x)\le 2$ for all $x\in V$.
\end{enumerate}
\end{proposition}

\begin{proof}
(i) $\Rightarrow$ (ii) is obvious by 
Proposition~\ref{03prop:cactus-like graph is a block graph},
and (ii) $\Rightarrow$ (i) by 
Lemma~\ref{01lem:H1 cap H2 <=1 for block graph}.

(ii) $\Rightarrow$ (iii). 
Suppose that $\mathrm{cd}(x)\ge 3$ for some $x\in V$.
Then there exist distinct maximal cliques $C_1,C_2$ and $C_3$
containing $x$. Since $C_1\cap C_2\cap C_3\neq\emptyset$,
we see that $C_1\sim C_2\sim C_3\sim C_1$ form a cycle of length $3$
in $\Gamma(G)$.
This contradicts that $\Gamma(G)$ is a tree.

(iii) $\Rightarrow$ (ii). 
Suppose that $\Gamma(G)$ is not a tree.
Then there exists a cycle in $\Gamma(G)$, say,
\[
C_1\sim C_2\sim \dotsb \sim C_s\sim C_1,
\]
where $s\ge3$ and $C_1, C_2, \dots, C_s$ are mutually distinct.
For $i=1,2,\dots,s$ we choose $x_i\in V$ such that
$C_i\cap C_{i+1}=\{x_i\}$, where we understand $C_{s+1}=C_1$.
We note that $x_1,x_2,\dots,x_s$ are mutually distinct.
In fact, if $x_i=x_j$ for some $1\le i<j\le s$,
we see that $x_i=x_j \in C_i\cap C_{i+1}\cap C_{j+1}$
and $\mathrm{cd}(x_i)\ge3$.
This contradicts the assumption of (iii).
Thus,  
\[
x_1\sim x_2\sim \dotsb \sim x_s\sim x_1,
\]
becomes a cycle in $G$ of length $s$.
We then see that $C_1\cup C_2\cup \dotsb \cup C_s$ is a 2-connected
subgraph.
Hence it is contained in a 2-connected component $H$,
which is a maximal clique by 
Lemma~\ref{01lem:2-connected comp = max clique}.
This contradicts that $C_1,\dots,C_s$ are distinct maximal cliques.
\end{proof}

\begin{remark}\normalfont
As an extremal case of condition (iii) in
Proposition~\ref{03prop:equivalent conditions for cactus-like},
we note that
a block graph $G=(V,E)$ is reduced to 
a finite or infinite complete graph if
$\mathrm{cd}(x)=1$ for all $x\in V$.
\end{remark}

\section{QEC of Block Graphs}
\label{04sec:QEC of Block Graphs}

\subsection{Upper Bound for QEC of Block Graphs}

We first present the key result required for our characterization of
cactus-like graphs.
The proof will be given in the next section.

\begin{theorem}\label{04thm:psi is CND}
Let $G=(V,E)$ be a block graph and 
assume that 
\[
p=\sup\{\mathrm{cd}(x)\,;\, x\in V\}<\infty.
\]
Define a matrix $\Psi=[\psi(x,y)]_{x,y\in V}$ by
\begin{equation}\label{04eqn:def of Psi on G}
\Psi=D-\frac{1}{p}(J-I),
\end{equation}
where $D=[d(x,y)]$ is the distance matrix of $G$,
$J$ the all-one matrix, and $I$ the identity matrix.
Or equivalently,
\begin{equation}\label{04eqn:def of psi=d-1/p}
\psi(x,y)=
\begin{cases}
0, & \text{if $x=y$}, \\
d(x,y)-\dfrac{1}{p}, & \text{if $x\ne y$},
\end{cases}
\qquad x,y\in V.
\end{equation}
Then $\Psi$ is conditionally negative definite.
\end{theorem}

The matrix $\Psi$ defined by \eqref{04eqn:def of Psi on G}
plays a crucial role in achieving our goal.
It is noteworthy that $\Psi$ has zeroes on the diagonal and
the off-diagonal entries are the distances shifted downwards by $1/p$.
As an immediate consequence of Theorem~\ref{04thm:psi is CND},
we derive an upper bound for $\mathrm{QEC}(G)$ for a
block graph $G$.

\begin{theorem}\label{04thm:upper bound QEC of block graph}
Let $G=(V,E)$ be a block graph with $|V|\ge2$ and 
assume that 
\[
p=\sup\{\mathrm{cd}(x)\,;\, x\in V\}<\infty.
\]
Then 
\begin{equation}\label{04eqn:upper bound of QEC(block graph)}
\mathrm{QEC}(G)\le -\frac{1}{p}\,.
\end{equation}
In particular, for any cactus-like graph $G$ we have 
$p\le 2$ and 
\begin{equation}\label{04eqn:upper bound of QEC(cactus-like graph)}
\mathrm{QEC}(G)\le -\frac{1}{2}\,.
\end{equation}
\end{theorem}

\begin{proof}
Take an arbitrary $f\in C_0(V)$ satisfying $\langle \bm{1},f\rangle=0$
and $\langle f,f\rangle=1$.
Using the matrix $\Psi$ introduced in Theorem~\ref{04thm:psi is CND},
we come to
\begin{align*}
\langle f, D f\rangle
&=\langle f, \Psi f\rangle +\frac{1}{p} \langle f, (J-I)f\rangle \\
&=\langle f, \Psi f\rangle +\frac{1}{p} 
 (\langle \bm{1},f\rangle^2 - \langle f, f\rangle) \\
&=\langle f, \Psi f\rangle -\frac{1}{p}.
\end{align*}
Moreover, since $\langle f, \Psi f\rangle \le 0$
by Theorem~\ref{04thm:psi is CND}, we obtain
\[
\langle f, D f\rangle
=\langle f, \Psi f\rangle -\frac{1}{p} 
\le -\frac{1}{p}\,,
\]
from which \eqref{04eqn:upper bound of QEC(block graph)} follows.
For a cactus-like graph $G$, we need only to note that
$\mathrm{cd}(x)\le 2$ for all $x\in V$, 
see Proposition~\ref{03prop:equivalent conditions for cactus-like}.
\end{proof}

\subsection{Upper Bound for QEC of Finite Block Graphs}

We will obtain a refinement of 
Theorem~\ref{04thm:upper bound QEC of block graph}
for finite block graphs
and complete the characterization of 
finite cactus-like graphs in terms of QEC.
We begin with the elementary fact that 
any finite block graph with at least two maximal cliques
contains an `end maximal clique.'

\begin{lemma}\label{04lem:end maximal clique}
Let $G=(V,E)$ be a finite block graph
and $\mathcal{C}$ the set of maximal cliques of $G$.
If $|\mathcal{C}|\ge2$, 
there exists $C\in \mathcal{C}$ such that
\begin{equation}\label{04eqn:end maximal clique (0)}
\bigg|\bigcup_{D\in\mathcal{C}\setminus\{C\}} D\cap C\bigg|=1.
\end{equation}
(In this case $C$ is called an end maximal clique.)
\end{lemma}

\begin{proof}
Since the clique graph $\Gamma(G)$ is connected
and $|\mathcal{C}|\ge2$ by assumption,
we see that for every $C\in \mathcal{C}$ there exists 
$D\in\mathcal{C}$ such that $D\sim C$, that is $|D\cap C|=1$.
Therefore, the left-hand side of
\eqref{04eqn:end maximal clique (0)} is at least one.
To prove the assertion by contradiction
we suppose that every $C\in \mathcal{C}$ satisfies
\begin{equation}\label{04eqn:end maximal clique}
\bigg|\bigcup_{D\in\mathcal{C}\setminus\{C\}} D\cap C\bigg|\ge2.
\end{equation}
Since $|C\cap D|\le 1$ for any $D\in\mathcal{C}\setminus\{C\}$,
\eqref{04eqn:end maximal clique} means that
there exist $D_1,D_2\in\mathcal{C}\setminus\{C\}$
such that $D_1\neq D_2$, $D_1\sim C\sim D_2$ 
and $D_1\cap C\neq D_2\cap C$.
Based on this we construct a walk in $\Gamma(G)$:  
\begin{equation}\label{04eqn:walk in Gamma(G)}
C_1\sim C_2 \sim C_3 \sim \dotsb \sim C_s \sim \dotsb,
\end{equation}
satisfying
\[
C_{s-1}\neq C_{s+1},
\qquad
C_{s-1}\cap C_s\neq C_s\cap C_{s+1},
\qquad s=2,3,\dots.
\]
In fact, we choose $C_1$ and $C_2$ in such a way that $C_1\sim C_2$,
and applying the aforementioned property 
we may find $C_3$ satisfying the desired condition.
Then we repeat this process.
However, since $G$ is a finite graph,
$\mathcal{C}$ is finite as a subset of the power set $2^{V}$.
Therefore, in \eqref{04eqn:walk in Gamma(G)} we may find
$s\ge4$ such that
$C_1, C_2, \dots, C_{s-1}$ are mutually distinct
and $C_s=C_t$ for some $1\le t\le s-3$.
Namely, we come to a cycle
\[
C_t\sim C_{t+1} \sim \dotsb \sim C_{s-1} \sim C_s=C_t
\]
of length $s-t\ge3$.
Then $\langle C_t\cup C_{t+1} \cup\dotsb \cup C_{s-1}\rangle$
becomes a 2-connected subgraph of $G$
and hence is contained in a 2-connected component.
But this 2-connected component is a clique
as $G$ is a block graph.
Thus the mutually distinct maximal cliques
$C_t,C_{t+1},\dots, C_{s-1}$ are contained in another
maximal clique, which is a contradiction.
\end{proof}

\begin{theorem}\label{04thm:upper bound for finite block graph}
Let $G$ be a finite block graph with two or more maximal cliques.
Then we have
\begin{equation}\label{04eqn:upper bound for finite block graph}
\mathrm{QEC}(G)<-\frac{1}{p}\,,
\qquad
p=\max\{\mathrm{cd}(x)\,;\, x\in V\}.
\end{equation}
\end{theorem}
 
\begin{proof}
Let $\mathcal{C}$ be the set of maximal cliques of $G$.
By assumption we have $|\mathcal{C}|\ge2$.
We prove \eqref{04eqn:upper bound for finite block graph}
by induction on the number $s=|\mathcal{C}|$.
For $s=2$ we have $G=K_m\star K_n$ for some $m,n\ge2$ and $p=2$.
Since 
\[
\mathrm{QEC}(K_m\star K_n)
=-\frac{1}{1+\sqrt{\Big(1-\dfrac{1}{m}\Big)\Big(1-\dfrac{1}{n}\Big)}}
<-\frac12\,,
\]
see Theorem~\ref{a:thm explicit value of Km star Kn} in Appendix,
the assertion \eqref{04eqn:upper bound for finite block graph}
certainly holds for $s=2$.
Let $s\ge3$ and assume that 
\eqref{04eqn:upper bound for finite block graph} holds up to $s-1$.
We need to show that \eqref{04eqn:upper bound for finite block graph}
holds for any finite block graph $G=(V,E)$ with $s$ maximal cliques.
We proceed by contradiction.

Suppose, for the sake of contradiction, that there exists 
a finite block graph $G$
with $s$ maximal cliques for which
\eqref{04eqn:upper bound for finite block graph} fails.
In light of Theorem~\ref{04thm:upper bound QEC of block graph}
we have 
\begin{equation}\label{03eqn:in proof QEC(G)}
\mathrm{QEC}(G)=-\frac1p.
\end{equation}
We note that $p=\max\{\mathrm{cd}(x)\,;\, x\in V\}\ge2$.
By Lemma~\ref{04lem:end maximal clique} there exists
an end maximal clique $W$ of $G$.
Let $a\in W$ be a unique vertex such that
$W\cap W_1=\{a\}$ for any $W_1\sim W$.
Set $V^\prime=V\setminus(W\setminus\{a\})$.
Then $G^\prime=\langle V^\prime\rangle
=G\setminus (W\setminus\{a\})$ becomes a
block graph with $s-1$ maximal cliques.
By the induction hypothesis we have
\begin{equation}\label{03eqn:in proof QEC(G')}
\mathrm{QEC}(G^\prime)<-\frac1p\,.
\end{equation}
On the other hand, by \eqref{03eqn:in proof QEC(G)}
we may choose $f\in C(V)$ such that
\begin{equation}\label{04eqn:in final proof (1)}
\langle f,Df\rangle=-\frac1p\,,
\qquad
\langle f,f\rangle=1\,,
\qquad
\langle \bm{1},f\rangle=0.
\end{equation}
Assume that $f(x)=0$ for any $x\in W\setminus\{a\}$.
Then $f$ is regarded as a function in $C(V^\prime)$.
Since $G^\prime$ is isometrically embedded in $G$,
the distance matrix of $G^\prime$ coincides 
with $D\!\restriction_{V^\prime}$.
Then \eqref{04eqn:in final proof (1)} implies
that $\mathrm{QEC}(G^\prime)\ge -1/p$,
which contradicts \eqref{03eqn:in proof QEC(G')}.
Therefore,
there exists some $x_0\in W\setminus\{a\}$ such that $f(x_0)\neq0$.

Taking $k\ge2$ we consider the star product 
$\Tilde{G}=(G,x_0)\star K_k$.
With the help of Proposition~\ref{02prop:Baskoro-Obata2021} 
and \eqref{03eqn:in proof QEC(G)} we come to
\begin{equation}\label{04eqn:in final proof (2)}
-\frac1p=\mathrm{QEC}(G)<\mathrm{QEC}(\Tilde{G}).
\end{equation}
On the other hand, since $\Tilde{G}=(\Tilde{V},\Tilde{E})$ is
a block graph satisfying $p=\max\{\mathrm{cd}(x)\,;\, x\in \Tilde{V}\}$,
by Theorem~\ref{04thm:upper bound QEC of block graph} we have
\begin{equation}\label{04eqn:in final proof (3)}
\mathrm{QEC}(\Tilde{G})\le -\frac1p.
\end{equation}
Since \eqref{04eqn:in final proof (2)}
and \eqref{04eqn:in final proof (3)} are incompatible,
we conclude that 
\eqref{04eqn:upper bound for finite block graph} holds
for any finite block graph $G$
with $s$ maximal cliques.
\end{proof}

\begin{remark}\label{04rem:cactus-like with p=1}\normalfont
If $G$ is a finite block graph with just one maximal clique,
it is a complete graph $K_n$ with $n\ge2$.
Then we have $p=1$ and $\mathrm{QEC}(G)=-1$ so that
\eqref{04eqn:upper bound for finite block graph} fails.
\end{remark}

\begin{remark}\normalfont
It is plausible that Theorem~\ref{04thm:upper bound for finite block graph}
remains valid for a block graph which is not necessarily finite
but has finitely many maximal cliques.
The detailed study will appear elsewhere.
\end{remark}

\subsection{Characterization of Cactus-like Graphs}

The main result (Theorem~\ref{01thm:main result})
of this paper is restated in the following

\begin{theorem}\label{04thm:main result}
Let $G=(V,E)$ be a finite connected graph with $|V|\ge2$.
Then the following conditions are equivalent:
\begin{enumerate}
\item[\upshape (i)] $G$ is cactus-like,
that is, every pair of distinct maximal cliques has 
at most one vertex in common and the clique graph $\Gamma(G)$ is a tree;
\item[\upshape (ii)]  $G$ is a block graph whose clique graph 
$\Gamma(G)$ is a tree;
\item[\upshape (iii)] $\mathrm{QEC}(G)<-1/2$.
\end{enumerate}
\end{theorem}

\begin{proof}
The equivalence (i) $\Leftrightarrow$ (ii) is already shown in
Proposition~\ref{03prop:equivalent conditions for cactus-like}.
As is mentioned in Introduction the implication
(iii) $\Rightarrow$ (i) was proved in \cite{Baskoro-Obata2024}.
Here we prove the implication (i) $\Rightarrow$ (iii).
Let $G=(V,E)$ be a cactus-like graph.
By Proposition~\ref{03prop:equivalent conditions for cactus-like} we have
$p=\max\{\mathrm{cd}(x)\,;\, x\in V\}\le2$.
If $p=2$, applying Theorem~\ref{04thm:upper bound for finite block graph},
we obtain $\mathrm{QEC}(G)< -1/2$.
If $p=1$, then $G$ is a complete graph and $\mathrm{QEC}(G)=-1<-1/2$,
see Remark \ref{04rem:cactus-like with p=1}.
\end{proof}

\subsection{Lower Bound for QEC of Block Graphs}

Let $G=(V,E)$ be a block graph with $|V|\ge2$.
Consider a vertex $x\in V$ and assume that $q=\mathrm{cd}(x)<\infty$.
Let $C_1,\dots,C_q$ be the maximal cliques containing $x$.
We see easily that
\[
\langle C_1, C_2,\dots,C_q\rangle
\cong K_{|C_1|}\star K_{|C_2|}\star\dotsb\star K_{|C_q|},
\]
which is an isometric subgraph of $G$.
Hence we have
\begin{equation}\label{04eqn:lower bound at x}
\mathrm{QEC}(G)\ge
\mathrm{QEC}(K_{|C_1|}\star K_{|C_2|}\star\dotsb\star K_{|C_q|}).
\end{equation}
The exact value of the right-hand side is described only implicitly
except some special cases,
see Theorem~\ref{a:thm explicit value of Km star Kn} in Appendix.
Neverthless, 
if we take the supremum of the right-hand side over $x\in V$,
we may obtain a better lower bound of $\mathrm{QEC}(G)$.
Here we derive a simple lower bound for $\mathrm{QEC}(G)$
from \eqref{04eqn:lower bound at x}.

\begin{proposition}\label{04prop:lower bound by n*}
Let $G=(V,E)$ be a block graph with $|V|\ge2$ and $x\in V$.
Assume that $q=\mathrm{cd}(x)<\infty$ and
let $C_1,\dots,C_q$ be the maximal cliques containing $x$.
If there exists a finite maximal clique containing $x$,
setting $n_*=\min\{|C_1|,|C_2|,\dots,|C_q|\}<\infty$,
we have
\begin{equation}\label{04eqn:lower bound by n*}
\mathrm{QEC}(G)\ge
-\frac{n_*}{(n_*-1)q+1}\,.
\end{equation}
In particular,
\begin{equation}\label{04eqn:lower bound by n*=2}
\mathrm{QEC}(G)\ge
-\frac{2}{q+1}\,.
\end{equation}
\end{proposition}

\begin{proof}
Since $K_{|C_1|}\star K_{|C_2|}\star\dotsb\star K_{|C_q|}$
contains the $q$-fold star power
$K_{n_*}^{\star q}$ as an isometric subgraph,
we see from \eqref{04eqn:lower bound at x} that
\[
\mathrm{QEC}(G)\ge
\mathrm{QEC}(K_{n_*}^{\star q})
=-\frac{n_*}{(n_*-1)q+1}\,,
\]
where the last equality owes to
Theorem~\ref{a:thm explicit value of star power} in Appendix.
Thus, we obtain \eqref{04eqn:lower bound by n*}.
Moreover, \eqref{04eqn:lower bound by n*=2} follows 
from \eqref{04eqn:lower bound by n*} by $n_*\ge2$.
Alternatively, \eqref{04eqn:lower bound by n*=2}
can be verified by observing that a complete bipartite graph $K_{1,q}$
is contained in $\langle C_1, C_2,\dots,C_q\rangle$ as
an isometric subgraph.
\end{proof}

\begin{proposition}\label{04prop:QEC(G)=-1/p}
Let $G=(V,E)$ be a block graph with $|V|\ge2$ and assume that
$p=\sup\{\mathrm{cd}(x)\,;\, x\in V\}<\infty$.
If there exists a vertex $x\in V$ such that
$\mathrm{cd}(x)=p$ and
every maximal clique containing $x$ is infinite,
then we have
\begin{equation}\label{04eqn:QEC(G)=-1/p}
\mathrm{QEC}(G)=-\frac{1}{p}\,.
\end{equation}
\end{proposition}

\begin{proof}
For any $n\ge2$,
$K_n^{\star p}$ is contained in
$K_\infty^{\star p}$ as an isometric subgraph.
We then see from \eqref{04eqn:lower bound at x} that
\[
\mathrm{QEC}(G)
\ge\mathrm{QEC}(K_n^{\star p})
=-\frac{n}{(n-1)p+1}
\]
for any $n\ge2$.
Letting $n\rightarrow\infty$ we obtain
\[
\mathrm{QEC}(G)\ge -\frac{1}{p}\,.
\]
With the help of Theorem~\ref{04thm:upper bound QEC of block graph},
we conclude that \eqref{04eqn:QEC(G)=-1/p} holds.
\end{proof}

\begin{proposition}
If a block graph $G=(V,E)$ satisfies
$\sup\{\mathrm{cd}(x)\,;\, x\in V\}=\infty$,
we have $\mathrm{QEC}(G)=0$.
\end{proposition}

\begin{proof}
For any $p\ge1$ there exists a vertex $x\in V$ such that
$\mathrm{cd}(x)\ge p$.
We may take maximal cliques $C_1,\dots,C_p$ containing $x$.
Since $\langle C_1,\dots,C_p \rangle$ contains 
$K_2^{\star p}\cong K_{1,p}$ as an isometric subgraph,
we have
\[
\mathrm{QEC}(G)
\ge \mathrm{QEC}(K_2^{\star p})
=-\frac{2}{p+1}\,.
\]
Letting $p\rightarrow\infty$, we obtain $\mathrm{QEC}(G)\ge0$.

On the other hand, Theorem~\ref{04thm:psi is CND} remains valid
understanding $\Psi=D$,
and so does Theorem~\ref{04thm:upper bound QEC of block graph}
understanding $\mathrm{QEC}(G)\le0$.
This is also justified by observing that
a block graph is obtained by repeated application of
star product of complete graphs.
Consequently, $\mathrm{QEC}(G)=0$.
\end{proof}

\section{Conditionally Negative Definite Kernels}
\label{05sec:Conditionally Negative Definite Kernels}

\subsection{Preliminaries}

Let $I$ be an arbitrary non-empty set.
Let $C(I;\mathbb{C})$ be the space of 
complex-valued functions $f:I\rightarrow \mathbb{C}$
and $C_0(I;\mathbb{C})$ the subspace of those with finite supports.
The inner product on $C_0(I)$ is defined by
\[
\langle f,g\rangle
=\sum_{i\in I} \overline{f(i)}\,g(i),
\qquad
f,g \in C_0(I;\mathbb{C}).
\]
Let $\bm{1}\in C(I;\mathbb{C})$ be the constant function defined by
$\bm{1}(i)=1$ for all $i\in I$.
Then, we write
\[
\langle \bm{1},f\rangle
=\sum_{i\in I}f(i),
\qquad
f \in C_0(I;\mathbb{C}).
\]
A complex-valued function $\psi:I\times I\rightarrow \mathbb{C}$
is called a \textit{kernel} on $I$.
In a natural manner, 
a function $f$ and a kernel $\psi$ on $I$ 
give rise to a column vector $f=[f(i)]$ and a matrix 
$\Psi=[\psi(i,j)]$, respectively.
Under this identification we use these notations freely.

A kernel $\psi$ on $I$ is called \textit{hermitian} if
\[
\psi(j,i)=\overline{\psi(i,j)},
\qquad i,j\in I.
\]
A kernel $\psi$ on $I$ is called \textit{positive definite} if
\[
\langle f, \Psi f\rangle
=\sum_{i,j\in I} \psi(i,j)\overline{f(i)}\,f(j)\ge0
\]
for any $f\in C_0(I;\mathbb{C})$.
Recall that a positive definite kernel is necessarily hermitian.
A kernel $\psi$ on $I$ is called
\textit{conditionally negative definite} if
it is hermitian and
\[
\langle f, \Psi f\rangle
=\sum_{i,j\in I} \psi(i,j)\overline{f(i)}\,f(j)\le0
\]
for any $f\in C_0(I;\mathbb{C})$ with $\langle\bm{1},f\rangle=0$.

\begin{lemma}\label{le:schwarz}
Let $\psi$ be a positive definite kernel on $I$.
Then we have
\[
\bigg|\sum_{i \in I} \psi(o,i) f(i)\bigg|^2
\le \psi(o,o)\langle f, \Psi f\rangle,
\qquad
o\in I, \quad f\in C_0(I;\mathbb{C}).
\]
\end{lemma}

\begin{proof}
We need only employ a standard argument using the Cauchy--Schwarz inequality,
based on the fact that
$\langle f+\lambda\delta_o, \Psi(f+\lambda\delta_o)\rangle\ge0$
for all $\lambda\in\mathbb{C}$,
where $\delta_o\in C_0(I)$ is defined by
$\delta_o(o)=1$ and $\delta_o(i)=0$ for $i\neq o$.
\end{proof}

\begin{lemma}\label{05lem:Schoenberg}
For a kernel $\psi$ on $I$ the following two conditions are equivalent:
\begin{enumerate}
\item[\upshape (i)] $\psi$ is conditionally negative definite;
\item[\upshape (ii)] the kernel defined 
by $e^{-t\psi}(i,j)=\exp(-t \psi(i,j))$ is positive definite 
for all $t>0$.
\end{enumerate}
\end{lemma}

The above result is well known in harmonic analysis
and is often referred to as Schoenberg's theorem
\cite{Berg-Christensen-Ressel1984,Parthasarathy-Schmidt1972}.
The kernel $e^{-t\psi}$ introduced in (ii) 
is also called the \textit{Hadamard exponential},
which coincides with the entry-wise exponential in a matrix form.

\subsection{Positive Definite Kernels on Block Graphs}
\label{04subsec:PD on block graphs}

Let $G=(V,E)$ be a block graph with $|V|\ge2$.
We will construct a positive definite kernel on $V$
by concatenating `local kernels' along a geodesic in $G$.
This method was first introduced for constructing
positive definite kernels on a tree \cite{Mlotkowski1992},
for relevant results, see also \cite{Mlotkowski1988,Mlotkowski1989,Valette1992}.

For a vertex $v_0\in V$ let $\mathcal{V}(v_0)$ denote 
the set of connected components of $G\setminus\{v_0\}$.
Then we come to a partition of $V$ defined by
\[
\mathcal{V}^\prime(v_0)=\mathcal{V}(v_0)\cup\{\{v_0\}\},
\]
and $V$ is decomposed into a disjoint union of non-empty subsets:
\begin{equation}\label{05eqn:partition of V}
V=\bigcup_{U\in\mathcal{V}^\prime(v_0)}U
 =\{v_0\}\cup \bigcup_{U\in\mathcal{V}(v_0)}U.
\end{equation}
Note that the number of connected components 
in $\mathcal{V}(v_0)$ coincides with $\mathrm{cd}(v_0)$.

Assume that each vertex $v_0\in V$ is assigned a
hermitian kernel on $\mathcal{V}^\prime(v_0)$.
We call the set $\{a_{v_0}\,;\, v_0\in V\}$
a \textit{system of local hermitian kernels}.
For a local hermitian kernel $a_{v_0}$ we define
$a(x,v_0,y)$ for $x,y\in V$ with
$d(x,v_0)\le 1$ and $d(y,v_0)\le 1$ in such a way that
\begin{equation}\label{05eqn:def of a}
a(x,v_0,y)=a_{v_0}(U,W),
\qquad U,W\in\mathcal{V}^\prime(v_0),
\quad x\in U, \quad y\in W,
\end{equation}
where $U$ and $W$ are
uniquely determined by \eqref{05eqn:partition of V}.
For \eqref{05eqn:def of a} we also write
\begin{equation}\label{eq:auvw}
a(x,v_0,y)=a(U,v_0,y)=a(x,v_0,W)=a(U,v_0,W)
\end{equation}
for convenience.

Let $\{a_{v_0}\,;\, v_0\in V\}$ be 
a system of local hermitian kernels,
and assume that $a_{v_0}$ is positive definite for all $v_0\in V$
and $a_{v_0}(U,U)=1$ for all $U\in\mathcal{V}^\prime(v_0)$.
Note the obvious relation:
\[
a(x,v_0,y)=1,
\qquad
x,y\in U,
\quad
U\in\mathcal{V}^\prime(v_0).
\]
We now define a `global' kernel on $V$
by concatenating local kernels along a geodesic.
Since $G$ is geodetic by Proposition~\ref{03prop:Block graph is geodetic},
for any pair of vertices $x,y\in V$ 
there exists a unique shortest path (the geodesic)
connecting $x$ and $y$,
denoted by
\begin{equation}\label{04eqn:x and y with distance n}
x=x_0\sim x_1 \sim x_2 \sim \dotsb \sim x_n=y,
\qquad n=d(x,y).
\end{equation}
Along this geodesic we define a kernel on $V$ by
\begin{equation}\label{04eqn:phi}
\phi(x,y)
=\prod_{k=0}^{n} a(x_{k-1},x_{k},x_{k+1})
\end{equation}
under the convention that $x_{-1}=x_0=x$ and $x_{n+1}=x_n=y$.
It is easy to see that $\phi$ becomes a hermitian kernel on $V$
satisfying $\phi(x,x)=1$ for all $x\in V$.
We then define an auxiliary kernel by
\begin{equation}\label{04eqn:beta}
\beta(x,y)
=\prod_{k=0}^{n-1} a(x_{k-1},x_{k},x_{k+1}),
\end{equation}
where \eqref{04eqn:x and y with distance n} is taken into account
with the convention that $x_{-1}=x_0=x$.
By definition $\beta(x,x)=1$ for all $x\in V$.
We see easily that
\[
\phi(x,y)=a(x,x,y)a(x,y,y)=\beta(x,y)\overline{\beta(y,x)}, 
\qquad x,y\in V, \quad x\sim y.
\]
Moreover, for $x,y\in V$ as in 
\eqref{04eqn:x and y with distance n} we have
\begin{equation}\label{04eqn:betaphibeta}
\phi(x,y)
=\beta(x,x_k)a(x_{k-1},x_k,x_{k+1})\overline{\beta(y,x_k)},
\qquad k=0,1,\dots,n.
\end{equation}
We mention two more formulas for $\phi$ in the following lemmas.

\begin{lemma}\label{04lem:beta-a-beta}
Let $v_0\in V$ and $U,W\in \mathcal{V}(v_0)$ with $U\neq W$.
Then we have
\[
\phi(x,y)=\beta(x,v_0)a_{v_0}(U,W)\overline{\beta(y,v_0)},
\qquad
x\in U,\quad y\in W.
\]
\end{lemma}

\begin{proof}
Noting that $x\neq v_0$ and $y\neq v_0$, we take the geodesics
\[
v_0=x_0\sim x_1\sim\dotsb\sim x_s=x,
\qquad
v_0=y_0 \sim y_1 \sim \dotsb\sim y_t=y.
\]
By assumption we have $s,t\ge1$, 
$x_1\in U$, $y_1 \in W$, and 
\[
x_s \sim x_{s-1} \sim \dotsb\sim x_1 \sim x_0=v_0=y_0
\sim y_1 \sim y_2 \sim\dotsb\sim y_t
\]
becomes the geodesic from $x$ to $y$.
Then the assertion follows by \eqref{04eqn:betaphibeta}.
\end{proof}

\begin{lemma}\label{04lem:beta-beta}
Let $v_0\in V$ and $U\in \mathcal{V}(v_0)$.
For any $x,y\in U$ satisfying $d(x,v_0)+d(y,v_0)=d(x,y)+1$ we have
\begin{equation}\label{eq:betabeta}
\phi(x,y)=\beta(x,v_0)\overline{\beta(y,v_0)}.
\end{equation}
\end{lemma}

\begin{proof}
Noting that $x\neq v_0$ and $y\neq v_0$, we take the geodesics
\[
v_0=x_0\sim x_1\sim\dotsb\sim x_s=x,
\qquad
v_0=y_0 \sim y_1 \sim \dotsb\sim y_t=y.
\]
By assumption we have $s,t\ge1$, $x_1\ne y_1$, $x_1\sim y_1$, and 
\[
x_s \sim x_{s-1} \sim \dotsb\sim x_1
\sim y_1 \sim y_2 \sim\dotsb\sim y_t
\]
becomes the geodesic from $x$ to $y$.
Then \eqref{eq:betabeta} follows by
\eqref{04eqn:phi} and \eqref{04eqn:beta}.
\end{proof}

We define a quadratic function by
\begin{equation}\label{04eqn:quadratic function associated to Phi}
Q(f)=\sum_{x,y\in V} \phi(x,y)f(x)\overline{f(y)},
\qquad f\in C_0(V).
\end{equation}
Note that $Q(f)=\langle \overline{f},\Phi \overline{f}\rangle$
in terms of the standard inner product
and the matrix notation $\Phi=[\phi(x,y)]$.
We then see that the kernel $\phi(x,y)$ is positive definite if 
$Q(f)\ge0$ for all $f\in C_0(V)$.
The use of $Q(f)$ is convenient for later argument.

To simplify the notation, we define a `localization' of $f\in C_0(V)$ by
\begin{equation}\label{04eqn:def of B-operation}
B_{v_0}f(U)
=\sum_{x\in U}\beta(x,v_0)f(x),
\qquad U\in\mathcal{V}^\prime(v_0),
\quad f\in C_0(V).
\end{equation}
Obviously, $B_{v_0}f\in C_0(\mathcal{V}^\prime(v_0))$.
We set
\begin{equation}\label{04eqn:local quadratic function}
Q_{v_0}(f)=\sum_{U,W\in\mathcal{V}^\prime(v_0)} 
 a_{v_0}(U,W) B_{v_0}f(U)\overline{B_{v_0}f(W)},
\qquad f\in C_0(V).
\end{equation}

\begin{lemma}\label{04lem:decomposition by components}
For any $v_0\in V$ and $f\in C_0(V)$ we have
\begin{equation}\label{04eqn:decomposition by components}
Q(f)-Q_{v_0}(f)
=\sum_{U\in\mathcal{V}(v_0)} \big\{Q(f|_U) -Q_{v_0}(f|_U)\big\},
\end{equation}
where $f|_U\in C_0(V)$ is defined by
$f|_U(x)=f(x)$ for $x\in U$ and $=0$ otherwise.
\end{lemma}

\begin{proof}
According to the partition of $V$ in \eqref{05eqn:partition of V}
the quadratic function 
\eqref{04eqn:quadratic function associated to Phi} becomes
\begin{align}
Q(f)=|f(v_0)|^2
&+\sum_{U\in\mathcal{V}(v_0)}\sum_{x,y\in U}\phi(x,y)f(x)\overline{f(y)}
\nonumber \\
&\qquad\quad +\sum_{\substack{U,W\in\mathcal{V}^\prime(v_0) \\ U \neq W}}
 \sum_{\substack{x\in U \\ y\in W}}\phi(x,y)f(x)\overline{f(y)}.
\label{04eqn:in proof Lemma 4.3 (1)}
\end{align}
The second term becomes
\begin{align}
&\sum_{U\in\mathcal{V}(v_0)}\sum_{x,y\in U}\phi(x,y)f(x)\overline{f(y)}
\nonumber \\
&\qquad
=\sum_{U\in\mathcal{V}(v_0)}\sum_{x,y\in V}\phi(x,y)f|_U(x)
 \overline{f|_U(y)}
=\sum_{U\in\mathcal{V}(v_0)} Q(f|_U).
\label{04eqn:in proof Lemma 4.3 (2)}
\end{align}
For the third term, in view of Lemma~\ref{04lem:beta-a-beta} 
and \eqref{04eqn:local quadratic function} we come to
\begin{align}
&\sum_{\substack{U,W\in\mathcal{V}^\prime(v_0) \\ U \neq W}}
\sum_{\substack{x\in U \\ y\in W}}\phi(x,y)f(x) \overline{f(y)}
\nonumber \\
&\qquad
 =\sum_{\substack{U,W\in\mathcal{V}^\prime(v_0) \\ U \neq W}}
 a_{v_0}(U,W) B_{v_0}f(U) \overline{B_{v_0}f(W)}
\nonumber \\
&\qquad = Q_{v_0}(f)
 -\sum_{U\in\mathcal{V}^\prime(v_0)} a_{v_0}(U,U) |B_{v_0}f(U)|^2.
\label{04eqn:in proof Lemma 4.3 (3)}
\end{align}
Since $B_{v_0}f|_U(W)=B_{v_0}f(U)$ for $W=U$ and $=0$ otherwise,
we see that
\begin{align}
&\sum_{U\in\mathcal{V}^\prime(v_0)} a_{v_0}(U,U) |B_{v_0}f(U)|^2
\nonumber \\
&\qquad = \sum_{U\in\mathcal{V}^\prime(v_0)}
 \sum_{U^\prime,W^\prime\in\mathcal{V}^\prime(v_0)} 
 a_{v_0}(U^\prime,W^\prime) B_{v_0}f|_U(U^\prime)\,
  \overline{B_{v_0}f|_U(W^\prime)}
\nonumber \\
&\qquad = |f(v_0)|^2
 + \sum_{U\in\mathcal{V}(v_0)} Q_{v_0}(f|_U).
 \label{04eqn:in proof Lemma 4.3 (4)}
\end{align}
Then \eqref{04eqn:decomposition by components} follows
by combining
\eqref{04eqn:in proof Lemma 4.3 (1)}--%
\eqref{04eqn:in proof Lemma 4.3 (4)}.
\end{proof}

Now we are ready to prove the main result of this section.

\begin{theorem}\label{th:pos}
Let $G=(V,E)$ be a block graph
and $\{a_{v_0}\,;\, v_0\in V\}$ a system of local hermitian kernels
satisfying $a_{v_0}(U,U)=1$ for any $U\in\mathcal{V}^\prime(v_0)$.
Define a hermitian kernel $\phi$ on $V$ as in \eqref{04eqn:phi}.
If $a_{v_0}$ is positive definite for all $v_0\in V$,
so is $\phi$.
\end{theorem}

\begin{proof}
We need to show that $Q(f)\ge0$ for any $f\in C_0(V)$,
where $Q(f)$ is a quadratic function defined in
\eqref{04eqn:quadratic function associated to Phi}.
Since $a_{v_0}$ is positive definite by assumption,
we see from \eqref{04eqn:local quadratic function} that
$Q_{v_0}(f)\ge0$ for any $v_0\in V$ and $f\in C_0(V)$.
We set
\begin{equation}\label{eq:delta}
\Delta(v_0,f)
=Q(f) - Q_{v_0}(f)
\end{equation}
Thus, for $Q(f)\ge0$
it is sufficient to show that $\Delta(v_0,f)\ge0$ 
for any $v_0\in V$ and $f\in C_0(V)$.
In fact, by induction we prove the following statement:
for each $n\ge0$ it holds that $\Delta(v_0,f)\ge0$ for any
$v_0\in V$ and $f\in C_0(V)$ satisfying
$\mathrm{supp}(f)\subset\{x\in V\,;\, d(v_0,x)\le n\}$.

For $n=0$, the statement is obvious
from
\[
\Delta(v_0,f)
=\phi(v_0,v_0)f(v_0)^2
 -a_{v_0}(\{v_0\},\{v_0\})f(v_0)^2=0.
\]
Let $n\ge1$ and assume that the statement is true up to $n-1$.
Now let $v_0\in V$ and $f\in C_0(V)$ with
$\mathrm{supp}(f)\subset \{x\in V\,;\, d(v_0,v)\le n\}$.
With the help of 
Lemma~\ref{04lem:decomposition by components} we have
\begin{equation}\label{eq:deltadecomp}
\Delta(v_0,f)=\sum_{U\in\mathcal{V}(v_0)} \Delta(v_0,f|_U).
\end{equation}
We will show that $\Delta(v_0,f|_{U_0})\ge0$
for any $U_0\in\mathcal{V}(v_0)$ by further decomposing
it into a sum of similar quadratic terms.
We set $K_0=\{v\in U_0:v_0\sim v\}$.
Then $K_0\cup\{v_0\}$ is a maximal clique of $G$,
and any geodesic connecting $x\in U_0$ and $v_0$
passes through a unique vertex in $K_0$.
We set
\[
V_v=\{x\in U_0:\, d(x,v_0)=d(x,v)+1\},
\qquad v\in K_0.
\]
In other words, $V_v$ is the set of vertices $x\in U_0$
such that the geodesic connecting $x$ and $v_0$ passes
through $v$, and thus $U_0$ is decomposed into a union
of disjoint subsets:
\[
U_0=\bigcup_{v\in K_0} V_v\,.
\]
For $v_1,v_2\in K_0$ with $v_1\ne v_2$ we have
\[
d(x,v_0)+d(y,v_0)=d(x,y)+1,
\qquad
x\in V_{v_1}, \quad y\in V_{v_2}.
\]
Then, by an argument similar to the proof of
Lemma~\ref{04lem:decomposition by components} with the help of
Lemma~\ref{04lem:beta-beta} we obtain
\begin{equation}\label{eq:deltav0}
\Delta(v_0,f|_{U_0})
=\sum_{v\in K_0}\Delta(v_0,f|_{V_v}),
\end{equation}
where
\begin{equation}\label{04eqn:in proof Theorem 4.4 (1)}
\Delta(v_0,f|_{V_v})
=\sum_{x,y\in V_v}\phi(x,y)f(x) \overline{f(y)}
  -\bigg|\sum_{x\in V_v}\beta(x,v_0)f(x)\bigg|^2.
\end{equation}

We next take $v_1\in K_0$ and consider $\Delta(v_0,f|_{V_{v_1}})$.
Set $\mathcal{V}_0(v_1)=\{X\in \mathcal{V}(v_1):v_0\notin X\}$.
For $X\in\mathcal{V}_0(v_1)$ we have
\[
\beta(x,v_0)=\beta(x,v_1) a(X,v_1,v_0),
\qquad x\in X.
\]
Then \eqref{04eqn:in proof Theorem 4.4 (1)} with $v$
being replaced by $v_1$ becomes
\begin{align}
\Delta(v_0,f|_{V_{v_1}})
&=\sum_{x,y\in V_{v_1}}\phi(x,y)f(x) \overline{f(y)}
-\bigg|\sum_{x\in V_{v_1}}\beta(x,v_0)f(x)\bigg|^2
\nonumber \\
&=\sum_{x,y\in V_{v_1}}\phi(x,y)f(x)\overline{f(y)}
 -\bigg|\sum_{X\in\mathcal{V}_0(v_1)}a(X,v_1,v_0)
  \sum_{x\in X}\beta(x,v_1)f(x)\bigg|^2
\nonumber \\
&=\sum_{x,y\in V_{v_1}}\phi(x,y)f(x)\overline{f(y)}
 -\bigg|\sum_{X\in\mathcal{V}_0(v_1)} a_{v_1}(X,\{v_0\})
  B_{v_1}f(X)\bigg|^2.
\label{04eqn:in proof Theorem 4.4 (2)}
\end{align}
Here we note that
\begin{equation}\label{04eqn:in proof Theorem 4.4 (3)}
\bigg|\sum_{X\in\mathcal{V}_0(v_1)}a_{v_1}(X,\{v_0\})B_{v_1}f(X)\bigg|^2
\le\sum_{X,Y\in\mathcal{V}_0(v_1)}
 a_{v_1}(X,Y)B_{v_1}f(X)\,\overline{B_{v_1}f(Y)},
\end{equation}
which follows from the Schwarz inequality (Lemma~\ref{le:schwarz})
and $a_{v_1}(\{v_0\},\{v_0\})=1$.
Applying \eqref{04eqn:in proof Theorem 4.4 (3)} to
\eqref{04eqn:in proof Theorem 4.4 (2)},
we obtain
\begin{align}
\Delta(v_0,f|_{V_{v_1}})
&\ge \sum_{x,y\in V_{v_1}}\phi(x,y)f(x) \overline{f(y)}
 -\sum_{X,Y\in\mathcal{V}_0(v_1)}
 a_{v_1}(X,Y)B_{v_1}f(X)\,\overline{B_{v_1}f(Y)}
\nonumber \\
&=\Delta(v_1,f|_{V_{v_1}})
\end{align}
Since $\Delta(v_1,f|_{V_{v_1}})\ge0$ by the induction hypothesis,
we have $\Delta(v_0,f|_{V_{v_1}})\ge0$,
that is, $\Delta(v_0,f|_{V_v})\ge0$ holds for any $v\in K_0$.
Then in view of \eqref{04eqn:in proof Theorem 4.4 (1)}
and \eqref{eq:deltav0} we have 
$\Delta(v_0,f|_{U_0})\ge0$ as desired.
\end{proof}

\subsection{Conditionally Negative Definite Kernels on Block Graphs}

Let $G=(V,E)$ be a block graph and
$\{c_{v_0}\,;\, v_0\in V\}$ a system of local hermitian kernels
satisfying $c_{v_0}(U,U)=0$ for any $U\in \mathcal{V}^\prime(v_0)$.
In a similar manner as in Section \ref{04subsec:PD on block graphs},
we define $c(x,v_0,y)$ for $x,y\in V$ with
$d(x,v_0)\le 1$ and $d(y,v_0)\le 1$,
see \eqref{05eqn:def of a}.
For any pair of vertices $x,y\in V$,
taking the geodesic
as in \eqref{04eqn:x and y with distance n} we set
\begin{equation}\label{05eqn:def of psi kernel}
\psi(x,y)=\sum_{i=0}^{n} c(x_{i-1},x_i,x_{i+1}).
\end{equation}
We see easily that $\psi$ becomes a hermitian kernel on $V$.

\begin{theorem}\label{th:neg}
Let $G=(V,E)$ be a block graph and
$\{c_{v_0}\,;\, v_0\in V\}$ a system of local hermitian kernels
satisfying $c_{v_0}(U,U)=0$ for any $U\in \mathcal{V}^\prime(v_0)$.
If $c_{v_0}$ is conditionally negative definite for all $v_0\in V$,
so is the hermitian kernel $\psi$
defined in \eqref{05eqn:def of psi kernel}.
\end{theorem}

\begin{proof}
For $v_0\in V$ and $t>0$ we define a kernel $a_{v_0,t}$
on $\mathcal{V}^\prime(v_0)$ by
\[
a_{v_0,t}(U,W)
=\exp(-t c_{v_0}(U,W)),
\qquad U,W\in \mathcal{V}^\prime(v_0).
\]
Since $c_{v_0}$ is conditionally negative definite,
the kernel $a_{v_0,t}$ is positive definite for all $t>0$
by Lemma~\ref{05lem:Schoenberg}.
It then follows from Theorem~\ref{th:pos} that
the hermitian kernel $\phi_t$ on $V$ defined by
\[
\phi_t(x,y)
=\prod_{i=0}^{n} a_t(x_{i-1},x_i,x_{i+1})
=\prod_{i=0}^{n} \exp(-t c(x_{i-1},x_i,x_{i+1}))
=\exp(-t \psi(x,y))
\]
is positive definite.
Again by Lemma~\ref{05lem:Schoenberg} we see that
$\psi$ is conditionally negative definite.
\end{proof}

\subsection{Proof of Theorem~\ref{04thm:psi is CND}}

We rely on the following elementary linear algebraic result.

\begin{lemma}\label{le:cnegativedef}
Let $m\ge1$.
For $u\in \mathbb{R}$ let $C_u$ be an $(m+1)\times (m+1)$ matrix
given in a block matrix form as
\[
C_u=\begin{bmatrix}
0 & u\bm{1}^T \\
u\bm{1} &  J-I 
\end{bmatrix},
\]
where $\bm{1}\in\mathbb{R}^m$ is the column vector 
whose entries are all one,
$J$ the $m\times m$ all-one matrix,
and $I$ the $m\times m$ identity matrix.
Then $C_u$ is conditionally negative definite if and only if 
$u\ge(m-1)/(2m)$.
\end{lemma}

\begin{proof}[Proof of Theorem~\ref{04thm:psi is CND}]
We define a system of local hermitian kernels $\{c_{v_0}\,;\, v_0\in V\}$ by
\begin{equation}\label{04eqn:def of c}
c_{v_0}(U,W)=c_{v_0}(W,U)=
\begin{cases}
0, & \text{if $U=W\in\mathcal{V}^\prime(v_0)$,} \\
\dfrac{p-1}{2p}, & \text{if $U=\{v_0\}$, $W\in\mathcal{V}(v_0)$}, \\
1, & \text{if $U,W\in\mathcal{V}(v_0)$, $U\ne W$.}
\end{cases}
\end{equation}
In fact, $c_{v_0}$ is a real symmetric kernel,
and $[c_{v_0}(U,W)]$ becomes an $(m+1)\times (m+1)$
matrix with $m=\mathrm{cd}(v_0)$.
Since $m\le p$, we see from Lemma~\ref{le:cnegativedef} that
the matrix $[c_{v_0}(U,W)]$ is conditionally negative definite.
Thus, $\{c_{v_0}\,;\, v_0\in V\}$ is a system of
local hermitian kernels such that
$c_{v_0}(U,U)=0$ for all $U\in \mathcal{V}^\prime(v_0)$
and $c_{v_0}$ is conditionally negative definite for all $v_0\in V$.

For $x,y\in V$ taking the geodesic
\[
x=x_0\sim x_1\sim \dots \sim x_n=y,
\qquad d(x,y)=n,
\]
we define
\[
\Tilde{\psi}(x,y)=\sum_{i=0}^{n} c(x_{i-1},x_i,x_{i+1}).
\]
It then follows from Theorem~\ref{th:neg} that
$\Tilde{\psi}$ is a conditionally negative definite kernel on $V$.
On the other hand, we see from \eqref{04eqn:def of c} that
\[
\Tilde{\psi}(x,x)=0,
\qquad x\in V,
\]
and if $d(x,y)=n\ge1$ we have 
\begin{align*}
\Tilde{\psi}(x,y)
&=c(x_0,x_0,x_1)+c(x_0,x_1,x_2)+
\dotsb+c(x_{n-2},x_{n-1},x_n)+c(x_{n-1},x_n,x_n) \\
&=\frac{p-1}{2p}+(n-1) +\frac{p-1}{2p}
=n-\frac{1}{p}
=d(x,y)-\frac{1}{p}\,.
\end{align*}
Therefore, $\Tilde{\psi}$ coincides with $\psi$ given in the
statement.
Consequently,
$\psi$ is a conditionally negative definite kernel on $V$.
\end{proof}

\section{Some Remarks}
\label{05sec:Some Remarks}

\subsection{Free Product of Complete Graphs}

\begin{figure}[hbt]
\begin{center}
\includegraphics[width=320pt]{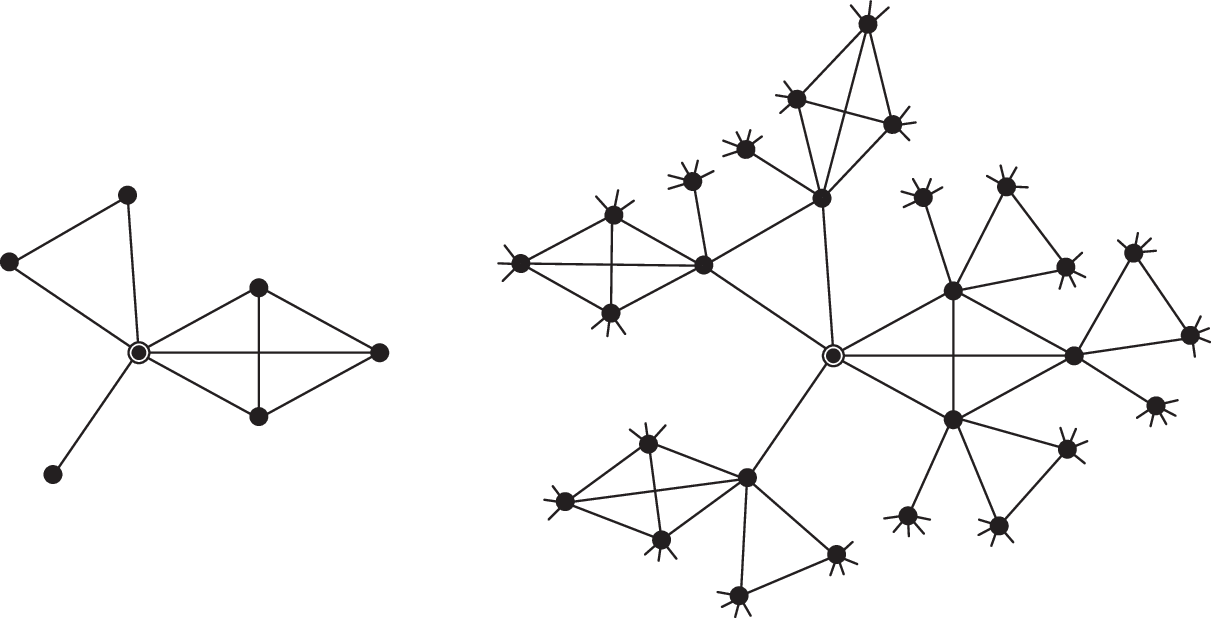}
\caption{Star product $K_2\star K_3 \star K_4$
and free product $K_2\circledast K_3 \circledast K_4$}
\label{06fig:forming free product graph}
\end{center}
\end{figure}

Let $p\ge2$.
For each $i=1,2,\dots,p$ let $G_i$ be
a finite or infinite complete garph
and consider their free product:
\begin{equation}\label{04eqn:free product graph}
G=G_1\circledast G_2 \circledast \cdots \circledast G_p\,.
\end{equation}
Let $n_i$ be the number of vertices of $G_i$ if it is finite,
otherwise we set $n_i=\infty$.
From now on, the free product \eqref{04eqn:free product graph}
is simply denoted by
\begin{equation}\label{04eqn:free product of complete graph}
G=K_{n_1}\circledast K_{n_2} \circledast \cdots \circledast K_{n_p}\,.
\end{equation}
An example is shown in Fig.~\ref{06fig:forming free product graph}.
The free product $G$ becomes an infinite block graph
which is homogeneous in the sense that
the subgraph induced by any vertex $x$ of $G$ and its
adjacent ones is isomorphic to the star product 
$K_{n_1}\star K_{n_2} \star \cdots \star K_{n_p}$.
We note that the degree and the clique degree are constant:
\[
\deg x=\sum_{i=1}^p(n_i-1),
\qquad
\mathrm{cd}(x)=p,
\qquad x\in V.
\]
In general, the free product is defined for
rooted graphs $(G_i,o_i)$ but for complete graphs $G_i$
we do not need to specify roots because of symmetry.
For a general definition and further properties of
free product graphs, see \cite{Carter-Tornier-Willis2020,Quenell1994,Znoiko1975}.

\begin{theorem}\label{04thm:QEC(G)<=-1/p}
Let $p\ge1$ and $n_i\ge2$ or $n_i=\infty$ for $i=1,2,\dots,p$.
If $n_i<\infty$ for some $i$,
setting $n_*=\min\{n_1,\dots,n_p\}<\infty$, we have
\begin{equation}\label{04eqn:qecfreeproduct}
-\frac{n_*}{p(n_*-1)+1}\le 
\mathrm{QEC}(K_{n_1} \circledast K_{n_2}
\circledast \dotsb \circledast K_{n_p})
\le -\frac{1}{p}.
\end{equation}
Moreover, if $n_i=\infty$ for all $i=1,2,\dots,p$,
we have
\begin{equation}\label{04eqn:p-fold power of infinite complete graphs}
\mathrm{QEC}(\,\overbrace{
K_\infty \circledast K_\infty
\circledast \dotsb \circledast K_\infty}^{\text{$p$ times}}\,)
=\mathrm{QEC}(K_\infty^{\circledast p})
=-\frac{1}{p}\,.
\end{equation}
\end{theorem}

\begin{proof}
The right half of \eqref{04eqn:qecfreeproduct} follows 
directly from Theorem~\ref{04thm:upper bound QEC of block graph},
and the left half from Proposition~\ref{04prop:lower bound by n*}.
For \eqref{04eqn:p-fold power of infinite complete graphs} 
it is sufficient to let $n_*\rightarrow\infty$
in \eqref{04eqn:qecfreeproduct}, 
see also Proposition~\ref{04prop:QEC(G)=-1/p}.
\end{proof}

\begin{theorem}\label{04thm:QEC=-1/2}
Let $n_i\ge2$ or $n_i=\infty$ for $i=1,2$.
Then the free product $K_{n_1} \circledast K_{n_2}$
is an infinite cactus-like graph and we have
\begin{equation}\label{04eqn:QEC=-1/2}
\mathrm{QEC}(K_{n_1} \circledast K_{n_2})=-\frac{1}{2}\,.
\end{equation}
\end{theorem}

\begin{figure}[hbt]
\begin{center}
\includegraphics[width=180pt]{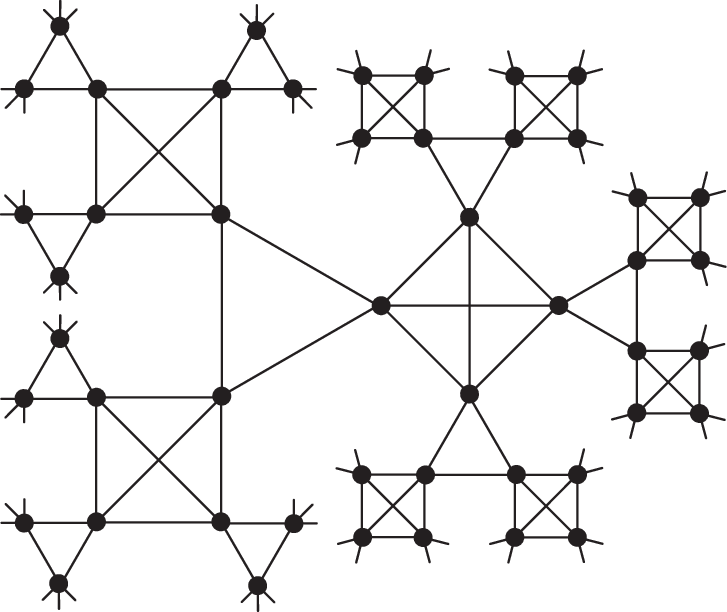}
\caption{Free product graph $K_3 \circledast K_4$}
\label{fig:free produc graph K3 free K4}
\end{center}
\end{figure}

\begin{proof}
Setting $p=2$ in Theorem~\ref{04thm:QEC(G)<=-1/p}
we obtain 
\begin{equation}\label{04eqn:in proof Thm 4.4 (0)}
\mathrm{QEC}(G)\le -\frac12.
\end{equation}
On the other hand, by construction
$G=K_{n_1} \circledast K_{n_2}$ contains 
$K_2 \circledast K_2\cong \mathbb{Z}$ (two-sided infinite path)
as an isometric subgraph.
Then, with the help of Propositions \ref{03prop:isometric subgraphs}
and \ref{03prop:QEC(Z)} we see that
\begin{equation}\label{04eqn:in proof Thm 4.4 (1)}
-\frac{1}{2}=\mathrm{QEC}(\mathbb{Z})
\le \mathrm{QEC}(G).
\end{equation}
Then the assertion follows from
\eqref{04eqn:in proof Thm 4.4 (0)}
and \eqref{04eqn:in proof Thm 4.4 (1)}.
An example of $K_{n_1} \circledast K_{n_2}$ is shown in 
Fig.~\ref{fig:free produc graph K3 free K4}.
\end{proof} 

\begin{remark}\normalfont
The free product graph
$G=K_{n_1}\circledast K_{n_2} \circledast \cdots \circledast K_{n_p}$
is well understood as the Cayley graph of the free product group
$\Gamma=\Gamma_1\circledast \Gamma_2 
\circledast \cdots \circledast \Gamma_p$,
where $\Gamma_i$ is a finite or infinite group 
consisting of $n_i$ elements,
with $\bigcup_{i=1}^p \Gamma_i\setminus\{e\}$ as the set of generators.
The length function on a free product group is 
important in harmonic analysis and has been studied considerably. 
Our argument for block graphs in Section 
\ref{05sec:Conditionally Negative Definite Kernels}
is an extension of the result on Cayley graphs of
free product groups developed in
\cite{Mlotkowski1988,Mlotkowski1989,Mlotkowski1992,Valette1992}.
\end{remark}

\subsection{Trees}

A graph is called a \textit{homogeneous tree}
of degree $p\ge2$, denoted by $T_p$,
if it is a tree and the degree of every vertex 
coincides with $p$.
For $p=1$ we may understand that $T_1=K_2$ though trivial.
Obviously, we have
\begin{equation}\label{05eqn:free product of K2}
T_p\cong K_2^{\circledast p} 
=K_2 \circledast K_2 \circledast \dotsb \circledast K_2
\quad\text{($p$ times)}.
\end{equation}
In particular, $T_2=K_2 \circledast K_2 \cong \mathbb{Z}$
becomes a two-sided infinite path.

\begin{proposition}\label{05prop:estimate of QEC(T)}
Let $T=(V,E)$ be a (finite or infinite) tree with $|V|\ge2$
and set $p=\sup\{\deg x\,;\,x\in V\}$.
If $1\le p<\infty$, we have
\begin{equation}\label{05eqn:estimate of trees}
-\frac{2}{p+1}\le \mathrm{QEC}(T)\le -\frac{1}{p}.
\end{equation}
If $p=\infty$, we have $\mathrm{QEC}(T)=0$.
\end{proposition}

\begin{proof}
Let $T=(V,E)$ be a tree with $|V|\ge2$ and assume that $1\le p<\infty$.
Since $T$ is isometrically embedded in 
the homogeneous tree $T_p$,
we see that $\mathrm{QEC}(T)\le \mathrm{QEC}(T_p)$.
Then the right half of \eqref{05eqn:estimate of trees} follows from
Theorem~\ref{04thm:QEC(G)<=-1/p}.
For the left half of \eqref{05eqn:estimate of trees} we take 
$x\in V$ with $\deg x=p$.
Let $y_1,\dots,y_p\in V$ be distinct vertices adjacent to $x$.
Since $T$ has no cycles, any pair of distinct $y_i$ and $y_j$ is not adjacent.
Then the induced subgraph $H=\langle x,y_1,\dots,y_p \rangle$
is isomorphic to the complete bipartite graph $K_{1,p}$.
Moreover, it is easy to see that $H$ is an isometric subgraph of $T$.
Thus, we obtain
\[
\mathrm{QEC}(T) \ge \mathrm{QEC}(H)
=\mathrm{QEC}(K_{1,p})=-\frac{2}{p+1}\,,
\]
where the explicit value of $\mathrm{QEC}(K_{1,p})$ is known
\cite{Obata-Zakiyyah2018}.
This proves the left half of \eqref{05eqn:estimate of trees}.
The rest of the assertion is obtained by letting $p\rightarrow\infty$
in \eqref{05eqn:estimate of trees}.
\end{proof}

\begin{remark}\normalfont
Let $T$ be a tree on $n\ge4$ vertices.
It is known \cite{Choudhury-Nandi2024} that
\[
\mathrm{QEC}(P_n) \le \mathrm{QEC}(T)
\le \mathrm{QEC}(K_{1, n-1})=-\frac{2}{n}\,.
\]
This is the best possible estimate in terms of the number of vertices 
in the sense that both equalities are attained.
Our result in Proposition~\ref{05prop:estimate of QEC(T)}
is an estimate of different type in terms of the maximal degree.
It is desirable to obtain a sharp estimate of
$\mathrm{QEC}(T)$ for a tree $T$ using other characteristics.
\end{remark}

\subsection{Further Remarks}

\begin{remark}\normalfont
Howorka \cite{Howorka1979} called a graph 
\textit{simple} if every pair of
distinct maximal cliques has at most one vertex in common.
As this terminology is apparently confusing today,
for the moment such a graph is called a H-graph.
$H$-graphs are interesting from the point of view of metric properties.
In this paper,
we have established the following inclusion 
for connected graphs:
\[
\{\text{cactus-like graphs}\}
\subset \{\text{block graphs}\}
\subset \{\text{geodetic graphs}\}
\subset \{\text{H-graphs}\}.
\]
The study of QEC along with metric properties would be 
an interesting future work.
\end{remark}

\begin{remark}\normalfont
The new term `cactus-like graph' is motivated by 
`cactus graph,' though the latter 
has been used with varying meanings in the literature.
According to the fairly standard definition,
a connected graph $G=(V,E)$ is called a \textit{cactus} if
every 2-conected component is an edge or a cycle.
In other words, a cactus is a connected graph in which
any two distinct cycles have at most one vertex in common,
or equivalently in which no edge lies on more than one cycle.
By definition the singleton graph $K_1$ is a cactus.
It is noted that cactus graphs and cactus-like graphs are distinct
though the two classes overlap.
When a cactus graph contains an even cycle,
it is clear that $\mathrm{QEC}(G)=0$.
When a cactus graph contains no even cycle,
further investigation is required.
\end{remark}

\setcounter{section}{1}
\section*{Appendix: Star Products of Complete Graphs}
\setcounter{equation}{0}
\setcounter{theorem}{0}
\renewcommand{\thesection}{\Alph{section}}
\renewcommand{\theequation}{\Alph{section}.\arabic{equation}}

Let $p\ge1$, $n_1\ge1,\dots,n_p\ge1$ and consider the star product:
\[
G=K_{n_1+1}\star K_{n_2+1}\star \dotsb \star K_{n_p+1}.
\]
The set of vertices of $G$ is given as the disjoint union:
\begin{equation}\label{a:disjoint union of vertices}
V=\{o\}\cup V_1\cup V_2\cup \dotsb \cup V_p\,,
\end{equation}
where $o$ is the root of $G$ and
$\{o\}\cup V_i$ is the vertex set of $K_{n_i+1}$.
Then,
according to the direct sum decomposition:
\begin{equation}\label{a:C(V) as direct sum}
C(V)=\mathbb{R}\oplus C(V_1)\oplus C(V_2) \oplus\dotsb \oplus C(V_p),
\end{equation}
any column vector $f\in C(V)$ is expressed as
\[
f=\begin{bmatrix}
f_0 \\ f_1 \\ \vdots \\ f_p
\end{bmatrix},
\qquad
f_0\in\mathbb{R},
\qquad f_i\in C(V_i)\cong \mathbb{R}^{n_i},
\quad i=1,2,\dots,p.
\]
Moreover, based on \eqref{a:disjoint union of vertices}
the distance matrix $D$ of $G$ is given as the block matrix
of the form:
\[
D=\begin{bmatrix}
0 & J & J & \cdots & J \\
J & J-I & 2J & \cdots & 2J \\
J & 2J & J-I & \cdots & 2J \\
\vdots & \vdots & \vdots & \ddots & \vdots \\
J & 2J & 2J & \cdots & J-I
\end{bmatrix},
\]
where $J$ is the all-one matrix and $I$ the identity matrix
of which the sizes are understood in the context.
By definition $\mathrm{QEC}(G)$ is the conditional maximum of 
$\langle f,Df\rangle$ 
subject to
\begin{equation}\label{a:constraints}
\langle f,f\rangle=f_0^2+\sum_{i=1}^p \langle f_i,f_i\rangle=1,
\qquad
\langle \bm{1},f\rangle=f_0+\sum_{i=1}^p \langle \bm{1},f_i\rangle=0.
\end{equation}
After simple algebra using \eqref{a:constraints} we come to
\begin{equation}\label{a:main quadratic function}
\langle f,Df\rangle
=- \sum_{i=1}^p
 (\langle f_i,f_i\rangle +\langle \bm{1},f_i\rangle^2)
\end{equation}
In this connection we recall the following

\begin{lemma}[\cite{Mlotkowski-Obata2020}]
\label{a:lemmma} 
Notations and assumptions being as above, 
let $a_i>0$ for $i=1,2,\dots, p$.
The conditional minimum of the quadratic function
\[
\sum_{i=1}^p a_i
 (\langle f_i,f_i\rangle +\langle \bm{1},f_i\rangle^2)
\]
subject to \eqref{a:constraints}
coincides with the minimal solution to
\[
\sum_{i=1}^p \frac{n_i}{a_in_i+a_i-\lambda}=\frac{1}{\lambda}\,.
\]
\end{lemma}

Then we have

\begin{theorem}\label{a:thm on star product}
Let $p\ge1$ and $n_1\ge1,\dots,n_p\ge1$.
Let $\lambda^*$ be the minimal solution to
\[
\sum_{i=1}^p \frac{n_i}{n_i+1-\lambda}=\frac{1}{\lambda}\,.
\]
Then we have
\[
\mathrm{QEC}(K_{n_1+1}\star K_{n_2+1}\star \dotsb \star K_{n_p+1})
=-\lambda^*\,.
\]
\end{theorem}

\begin{proof}
By direct application of Lemma~\ref{a:lemmma}
to \eqref{a:main quadratic function}.
\end{proof}

\begin{theorem}\label{a:thm explicit value of Km star Kn}
For $m\ge2$ and $n\ge2$ we have
\[
\mathrm{QEC}(K_m\star K_n)
=-\frac{mn}{mn+\sqrt{mn(m-1)(n-1)}}\,.
\]
\end{theorem}

\begin{proof}
Straightforward by
setting $n_1=m-1$ and $n_2=n-1$ in Theorem~\ref{a:thm on star product}.
\end{proof}

\begin{theorem}\label{a:thm explicit value of star power}
Let $p\ge1$ and $n\ge2$.
Then for the $p$-fold star power
$K_n^{\star p}=K_n\star \dotsb\star K_n$
($p$ times) we have
\[
\mathrm{QEC}(K_n^{\star p})=-\frac{n}{p(n-1)+1}\,.
\]
\end{theorem}

\begin{proof}
Straightforward by
setting $n_i=n-1$ in Theorem~\ref{a:thm on star product}.
\end{proof}

\end{document}